\documentclass[dspr]{apjrnl}
\usepackage{amssymb}
\usepackage{amsmath}
\usepackage{amsfonts}
\usepackage[doublespacing]{setspace}

{\theoremstyle{plain}

{\theoremstyle{remark}

}
{\theoremstyle{definition}

}

\iftwocol\AtEndDocument{\end{multicols}}\fi
\input{tcilatex}
\begin{document}

\title{Inverse Problems for Ultrahyperbolic Schr\"{o}dinger Equations}
\author{Fikret G\"{o}lgeleyen and \"{O}zlem Kaytmaz \\
Department of Mathematics, Bulent Ecevit University, Zonguldak 67100 Turkey 
\\
E-mail: f.golgeleyen@beun.edu.tr, ozlem.kaytmaz@beun.edu.tr}
\maketitle

\begin{abstract}
\textbf{Abstract.} In this paper, we establish a global Carleman estimate
for an Ultrahyperbolic Schr\"{o}dinger equation. Moreover, we prove H\"{o}%
lder stability for the inverse problem of determining a coefficient or a
source term in the Ultrahyperbolic Schr\"{o}dinger equation by some lateral
boundary data.

\end{abstract}

\section{\textbf{Introduction}}

Let $n,m\in 
\mathbb{N}
,$ $T>0$ and let $D\subset 
\mathbb{R}
^{n}$ be a bounded domain with smooth boundary $\partial D$ and $G=\left\{
y\in 
\mathbb{R}
^{m};\text{ }\left\vert y\right\vert <L\right\} $ for $L>0.$

We set $Q=D\times G\times (0,T)$, $\Sigma =\partial D\times G\times (0,T)$
and $i=\sqrt{-1}.$

We consider the ultrahyperbolic\textbf{\ }Schr\"{o}dinger equation%
\begin{equation}
i\partial _{t}v(x,y,t)+\Delta _{y}v(x,y,t)-\Delta
_{x}v(x,y,t)-p(x,y)v(x,y,t)=0,\text{ }(x,y,t)\in Q,  \tag{1.1}  \label{1.1}
\end{equation}%
with the following initial and Dirichlet boundary data%
\begin{eqnarray}
v(x,y,0) &=&a(x,y),\text{ }(x,y)\in D\times G,  \TCItag{1.2}  \label{1.2} \\
v(x,y,t) &=&0,\text{ }(x,y,t)\in \Sigma .  \TCItag{1.3}  \label{1.3}
\end{eqnarray}%
Throughout this paper, we use the following notations:%
\begin{eqnarray*}
\partial _{t} &=&\frac{\partial }{\partial t},\text{ }\partial _{x_{i}}=%
\frac{\partial }{\partial x_{i}},\text{ }\partial _{y_{j}}=\frac{\partial }{%
\partial y_{j}},\text{ }\nabla _{x}=\left( \partial _{x_{1}},\partial
_{x_{2}},...,\partial _{x_{n}}\right) ,\text{ } \\
\nabla _{y} &=&\left( \partial _{y_{1}},\partial _{y_{2}},...,\partial
_{y_{m}}\right) ,\text{ }\Delta _{x}=\sum\limits_{i=1}^{n}\partial
_{x_{i}}^{2},\text{ }\Delta _{y}=\sum\limits_{j=1}^{m}\partial _{y_{j}}^{2}.
\end{eqnarray*}

Let $v=v(p)$ satisfy (\ref{1.1})--(\ref{1.3}). We discuss the following
coefficient inverse problem.

\textbf{Problem 1 }Determine the coefficient $p(x,y),$ $\left( x,y\right)
\in D\times G$ in (\ref{1.1})--(\ref{1.3}) by the extra data $\left.
\partial _{\nu }v(p)\right\vert _{\Sigma }$, where $\nu \in 
\mathbb{R}
^{n}$ denotes the unit outward normal vector to $\partial D$ and $\partial
_{\nu }v=\nabla v\cdot \nu $ is the normal derivative.

Ultrahyperbolic Schr\"{o}dinger equations arise in several applications, for
example in water wave problems, [10, 11, 14, 34, 35] and in higher
dimensions as completely integrable models, see [1, 25]. There have been
limited number of studies on the direct problems for these equations. The
local well posedness of the initial value problem was investigated in
[19-21]. To our best knowledge there is no result available in the
mathematical literature related to the inverse problems for ultrahyperbolic
Schr\"{o}dinger equations. In this work,\ we obtain a Carleman estimate and
prove conditional H\"{o}lder stability for the inverse problem of
determining a coefficient or a source term in ultrahyperbolic Schr\"{o}%
dinger equation.

T. Carleman [8] established the first Carleman estimate in 1939 for proving
the unique continuation for a two-dimensional elliptic equation. In 1954, C.
M\"{u}ller extended Carleman's result to $%
\mathbb{R}
^{n},$ [18]. After that A. P. Calder\'{o}n [7] and L. H\"{o}rmander [13]
improved these results based on the concept of pseudo-convexity.

In the theory of inverse problems, Carleman estimates were firstly
introduced by A. L. Bukgeim and M. V. Klibanov in [6]. After that, there
have been many works relying on that method with modified arguments. Puel
and Yamamoto [28], Isakov and Yamamoto [17], Imanuvilov and Yamamoto [15,
16], Bellassoued and Yamamoto [5], Klibanov and Yamamoto [24] have obtained
various stability estimates for inverse problems for hyperbolic equation.

We refer to Yamamoto [31] for a comprehensive survey about the stability and
observability results for inverse problems for parabolic equations.

Inverse problems for ultrahyperbolic equations were considered in [2, 26,
29], where the unique continuation and stability were proved by using the
Carleman estimates. G\"{o}lgeleyen and Yamamoto [12] proved conditional H%
\"{o}lder stability for some inverse problems for ultrahyperbolic equation.

If $n=0$ then (\ref{1.1}) is a classical Schr\"{o}dinger equation which
describes the evalution of wave function of a charged particle under the
influence of electrical potential $p.$ As for the classical Schr\"{o}dinger
equation, Baudouin and Puel in [3] established a global Carleman estimate
and proved the uniqueness and Lipshchitz stability based on the idea by
Imanuvilov and Yamamoto [16]. This result was improved by Mercado et al.
[27] under a relaxed pseudoconvexity condition. In [3, 27], the main
assumption is that the part of the boundary where the measurement is made
satisfies a geometric condition related to geometric optics condition for
the observability. This geometric condition was relaxed in Bellassoued and
Choulli [4] under the assumption that the potential is known in a
neigborhood of the boundary of the spatial domain. Yuan and Yamamoto [33]
obtained a Carleman estimate with a regular weight function in Sobolev
spaces of negative orders. They proved the uniqueness in the inverse problem
of determining $L^{p}$ potentials and obtained an $L^{2}$ level
observability inequality and unique continuation results for the Schr\"{o}%
dinger equation. Cristofol and Soccorsi [9] considered the inverse problem
of determining time-dependent coefficient of the Schr\"{o}dinger equation
from a finite number of Neumann data. Kian et al. [22] extended the
stability result of [3] to the case of unbounded domains.

This paper consists of four sections. The rest of the paper is organized as
follows. In Section 2, the main result of this paper (Theorem 1) is
presented. In Section 3, a Carleman estimate (Proposition 1) which will be
used in the proof of our main result is established. Finally, Section 4 is
devoted to the proof of Theorem 1.

\section{\textbf{Main Result}}

Problem 1 can be reduced to an inverse source problem. For this aim, let $%
v(p)$ and $v(q)$ be two solutions of (\ref{1.1})--(\ref{1.3}) with the
coefficients $p$ and $q$ respectively.

Then the difference $u=v(p)-v(q)$ satisfies%
\begin{eqnarray}
Au &=&i\partial _{t}u(x,y,t)+\Delta _{y}u(x,y,t)-\Delta
_{x}u(x,y,t)-p(x,y)u(x,y,t)  \notag \\
&=&f(x,y)R(x,y,t),\text{ }(x,y,t)\in Q,  \TCItag{2.1}  \label{2.1} \\
u(x,y,0) &=&0,\text{ }(x,y)\in D\times G,  \TCItag{2.2}  \label{2.2} \\
u(x,y,t) &=&0,\text{ }(x,y,t)\in \Sigma  \TCItag{2.3}  \label{2.3}
\end{eqnarray}%
with $f(x,y)=p(x,y)-q(x,y)$ and $R=v(q)(x,y,t).$

We consider the following inverse source problem:

\textbf{Problem 2 }Let $p,R$ be given suitably. Then determine $f(x,y),$ $%
\left( x,y\right) \in D\times G$ by the extra data $\left. \partial _{\nu
}u\right\vert _{\Sigma }.$

Here we do not assume the uniqueness of $v(p)$ and $v(q)$ but their
existence. We have the following result for Problem 2.

\textbf{Theorem} \textbf{1} Let $p\in L^{\infty }\left( D\times
\{|y|<2L\}\right) $ and $u$ satisfy (\ref{2.1})--(\ref{2.3}) in $D\times
\{|y|<2L\}\times (-T,T).$

We assume that%
\begin{equation*}
u=0\text{ on }\partial D\times \{|y|<2L\}\times (-T,T),
\end{equation*}%
\begin{equation*}
\left\Vert \partial _{t}^{k}u\right\Vert _{H^{2}\left( D\times
\{|y|<2L\}\times (-T,T)\right) }\leq M,\text{ }k=1,2\text{ and}
\end{equation*}%
\begin{equation*}
R(x,y,0)\in 
\mathbb{R}
\text{ or }iR(x,y,0)\in 
\mathbb{R}
\text{ a. e. in }\left( D\times \{|y|<2L\}\right) ,
\end{equation*}%
\begin{equation*}
R\in H^{2}\left( -T,T;L^{\infty }\left( D\times \{|y|<2L\}\right) \right) ,
\end{equation*}%
\begin{equation*}
\left\Vert \partial _{t}^{k}R\right\Vert _{L^{2}(-T,T;L^{\infty }(D\times
\{|y|<2L\}))}\leq M,\text{ }k=1,2.
\end{equation*}%
We further assume that there exists a constant $r_{0}>0$ such that%
\begin{equation*}
|R(x,y,0)|\geq r_{0},\text{ }x\in \overline{D},\text{ }|y|\leq 2L
\end{equation*}%
and $\alpha >0$ is sufficiently small and%
\begin{equation*}
L>\frac{1}{\sqrt{\alpha }}\max_{x\in \overline{D}}\left\vert
x-x_{0}\right\vert .
\end{equation*}%
Then for any small $\epsilon >0$, for all real valued $f\in L^{2}(D\times
\{|y|<L\})$ there exist constants $C>0$ and $\theta \in (0,1)$ depending on $%
\epsilon ,$ $M,$ $x_{0},$ such that%
\begin{equation*}
\left\Vert f\right\Vert _{_{L^{2}(D\times \{|y|<L-\epsilon \})}}\leq
C\sum\limits_{k=1}^{2}\left\Vert \partial _{\nu }\partial
_{t}^{k}u\right\Vert _{_{L^{2}\left( \partial D_{+}\times \{|y|<2L\}\times
(-T,T)\right) }}^{\theta }.
\end{equation*}

\section{\textbf{Key Carleman Estimate}}

In this section, we show a Carleman estimate for the ultrahyperbolic\textbf{%
\ }Schr\"{o}dinger equation which will be used in the proof of our main
result.

Let $\Omega =D\times G\times (-T,T),$ $\Gamma _{x}=\partial D\times G\times
(-T,T)$, $\Gamma _{y}=D\times \partial G\times (-T,T)$ and let $\partial
\Omega =\Gamma _{x}\cup \Gamma _{y}.$

Let $x_{0}\notin \overline{D},$ $y_{0}\in 
\mathbb{R}
^{m}$ and $\alpha ,$ $\beta \in \left( 0,1\right) ,$ we set the weight
function%
\begin{equation}
\varphi \left( x,y,t\right) =e^{\gamma \psi \left( x,y,t\right) },  \tag{3.1}
\label{3.1}
\end{equation}%
where 
\begin{equation}
\psi \left( x,y,t\right) =|x-x_{0}|^{2}-\alpha |y-y_{0}|^{2}-\beta
\left\vert t\right\vert ^{2},  \tag{3.2}  \label{3.2}
\end{equation}%
$\gamma >0$ is a parameter. Moreover, we set the geometric condition%
\begin{equation}
\partial D_{+}=\left\{ x\in \partial D;\text{ }(x-x_{0})\cdot \nu \geq
0\right\}  \tag{3.3}  \label{3.3}
\end{equation}%
for $x_{0}\notin \overline{D}.$

We set%
\begin{gather}
Lu:=i\partial _{t}u\left( x,y,t\right) +\Delta _{y}u\left( x,y,t\right)
-\Delta _{x}u\left( x,y,t\right) +\sum\limits_{i=1}^{n}a_{i}\left(
x,y,t\right) u_{x_{i}}  \notag \\
+\sum\limits_{j=1}^{m}b_{j}\left( x,y,t\right) u_{y_{j}}+a_{0}\left(
x,y,t\right) u,  \tag{3.4}  \label{3.4}
\end{gather}%
\bigskip where $a_{i},$ $b_{j}\in L^{\infty }(\Omega ),$ $0\leq i\leq n,$ $%
1\leq j\leq m.$

\textbf{Proposition} \textbf{1} Let us assume that $0<\alpha ,$ $\beta <1$
be small and $\gamma >0$ be sufficiently large, and let%
\begin{equation}
|x-x_{0}|^{2}-\alpha ^{2}|y|^{2}-\beta ^{2}\left\vert t\right\vert
^{2}>\delta _{0}^{2},\text{ }\left( x,y,t\right) \in \Omega  \tag{3.5}
\label{3.5}
\end{equation}%
with some $\delta _{0}>0.$ Then there exist constants $C>0$ and $s_{0}>0$
such that%
\begin{eqnarray*}
&&\int_{\Omega }(s\left\vert \nabla _{y}u\right\vert ^{2}+s\left\vert \nabla
_{x}u\right\vert ^{2}+s^{3}\left\vert u\right\vert ^{2})e^{2s\varphi }dxdydt
\\
&\leq &C\int_{\Omega }\left\vert Lu\right\vert ^{2}e^{2s\varphi
}dxdydt+C\int_{\partial D_{+}\times G\times (-T,T)}s\left\vert \partial
_{\nu }u\right\vert ^{2}e^{2s\varphi }dS_{x}dydt
\end{eqnarray*}%
for all $s\geq s_{0},$ provided that 
\begin{eqnarray}
Lu &\in &L^{2}\left( \Omega \right) ,\text{ }u\in H^{2}\left( \Omega \right)
,  \notag \\
u\left( x,y,t\right) &=&0,\text{ }\left( x,y,t\right) \in \Gamma _{x}, 
\notag \\
u\left( x,y,t\right) &=&\left\vert \nabla _{y}u\left( x,y,t\right)
\right\vert =0,\text{ }\left( x,y,t\right) \in \Gamma _{y},  \notag \\
u\left( x,y,T\right) &=&u\left( x,y,-T\right) =0,\text{ }\left( x,y\right)
\in D\times G.  \TCItag{3.6}  \label{3.6}
\end{eqnarray}

\textbf{Proof }Let us set%
\begin{equation}
L_{0}u:=i\partial _{t}u+\Delta _{y}u-\Delta _{x}u=F,  \tag{3.7}  \label{3.7}
\end{equation}%
and%
\begin{equation}
z\left( x,y,t\right) =e^{s\varphi }u\left( x,y,t\right) ,\ P_{s}z\left(
x,y,t\right) =e^{s\varphi }L_{0}u.  \tag{3.8}  \label{3.8}
\end{equation}%
By (\ref{3.8}), we calculate%
\begin{eqnarray*}
P_{s}z &=&e^{s\varphi }L_{0}u \\
&=&i\partial _{t}z-is\partial _{t}\varphi z+\Delta _{y}z-\Delta
_{x}z-2s\left( \nabla _{y}\varphi \cdot \nabla _{y}z-\nabla _{x}\varphi
\cdot \nabla _{x}z\right) \\
&&-s\left( \Delta _{y}\varphi -\Delta _{x}\varphi \right) z+s^{2}\left(
\left\vert \nabla _{y}\varphi \right\vert ^{2}-\left\vert \nabla _{x}\varphi
\right\vert ^{2}\right) z.
\end{eqnarray*}%
Then we have%
\begin{equation}
P_{s}z+is\partial _{t}\varphi z=P_{s}^{+}z+P_{s}^{-}z,  \tag{3.9}
\label{3.9}
\end{equation}%
where%
\begin{eqnarray}
P_{s}^{+}z &=&i\partial _{t}z+\Delta _{y}z-\Delta _{x}z+s^{2}\left(
\left\vert \nabla _{y}\varphi \right\vert ^{2}-\left\vert \nabla _{x}\varphi
\right\vert ^{2}\right) z,  \TCItag{3.10}  \label{3.10} \\
P_{s}^{-}z &=&-2s\left( \nabla _{y}\varphi \cdot \nabla _{y}z-\nabla
_{x}\varphi \cdot \nabla _{x}z\right) -s\left( \Delta _{y}\varphi -\Delta
_{x}\varphi \right) z,  \TCItag{3.11}  \label{3.11}
\end{eqnarray}%
with the conventions $z\cdot z^{^{\prime
}}=\sum\limits_{i=1}^{N}z_{i}z_{i}^{^{\prime }}$ for all $z=\left(
z_{1},...,z_{N}\right) \in 
\mathbb{C}
^{N},$ $z^{^{\prime }}=\left( z_{1}^{^{\prime }},...,z_{N}^{^{\prime
}}\right) \in 
\mathbb{C}
^{N}.$

Then we have%
\begin{equation}
\left\Vert P_{s}z+is\partial _{t}\varphi z\right\Vert _{L^{2}\left( \Omega
\right) }^{2}=\left\Vert P_{s}^{+}z\right\Vert _{L^{2}\left( \Omega \right)
}^{2}+\left\Vert P_{s}^{-}z\right\Vert _{L^{2}\left( \Omega \right) }^{2}+2%
\func{Re}\left( P_{s}^{+}z,P_{s}^{-}z\right) _{L^{2}\left( \Omega \right) },
\tag{3.12}  \label{3.12}
\end{equation}%
where $\func{Re}\left( z\right) $ is the real part of $z.$\ Now, we
calculate the last term in (\ref{3.12}) by using (\ref{3.10}) and (\ref{3.11}%
), we obtain%
\begin{equation}
2\func{Re}\left( P_{s}^{+}z,P_{s}^{-}z\right) _{L^{2}\left( \Omega \right)
}=I_{1}+I_{2}+I_{3}+I_{4}+I_{5}+I_{6}+I_{7}+I_{8},  \tag{3.13}  \label{3.13}
\end{equation}%
where%
\begin{eqnarray*}
I_{1} &=&-4s\func{Re}\int_{\Omega }i\partial _{t}z\left( \nabla _{y}\varphi
\cdot \nabla _{y}\overline{z}-\nabla _{x}\varphi \cdot \nabla _{x}\overline{z%
}\right) dxdydt, \\
I_{2} &=&-2s\func{Re}\int_{\Omega }i\partial _{t}z\left( \Delta _{y}\varphi
-\Delta _{x}\varphi \right) \overline{z}dxdydt, \\
I_{3} &=&-4s\func{Re}\int_{\Omega }\Delta _{y}z\left( \nabla _{y}\varphi
\cdot \nabla _{y}\overline{z}-\nabla _{x}\varphi \cdot \nabla _{x}\overline{z%
}\right) dxdydt, \\
I_{4} &=&-2s\func{Re}\int_{\Omega }\Delta _{y}z\left( \Delta _{y}\varphi
-\Delta _{x}\varphi \right) \overline{z}dxdydt, \\
I_{5} &=&4s\func{Re}\int_{\Omega }\Delta _{x}z\left( \nabla _{y}\varphi
\cdot \nabla _{y}\overline{z}-\nabla _{x}\varphi \cdot \nabla _{x}\overline{z%
}\right) dxdydt, \\
I_{6} &=&2s\func{Re}\int_{\Omega }\Delta _{x}z\left( \Delta _{y}\varphi
-\Delta _{x}\varphi \right) \overline{z}dxdydt, \\
I_{7} &=&-4s^{3}\func{Re}\int_{\Omega }\left( \left\vert \nabla _{y}\varphi
\right\vert ^{2}-\left\vert \nabla _{x}\varphi \right\vert ^{2}\right)
z\left( \nabla _{y}\varphi \cdot \nabla _{y}\overline{z}-\nabla _{x}\varphi
\cdot \nabla _{x}\overline{z}\right) dxdydt, \\
I_{8} &=&-2s^{3}\func{Re}\int_{\Omega }\left( \left\vert \nabla _{y}\varphi
\right\vert ^{2}-\left\vert \nabla _{x}\varphi \right\vert ^{2}\right)
z\left( \Delta _{y}\varphi -\Delta _{x}\varphi \right) \overline{z}dxdydt,
\end{eqnarray*}%
and $\overline{z}$ is the conjugate of $z$.

Now, we shall estimate the terms $I_{k}$, $1\leq k\leq 8,$ using the
integration by parts and the condition $z\left( x,y,\pm T\right) =0.$ Then
we have%
\begin{eqnarray}
I_{1} &=&-4s\func{Re}\int\nolimits_{\Omega }i\partial _{t}z\nabla
_{y}\varphi \cdot \nabla _{y}\overline{z}dxdydt+4s\func{Re}%
\int\nolimits_{\Omega }i\partial _{t}z\nabla _{x}\varphi \cdot \nabla _{x}%
\overline{z}dxdydt  \notag \\
&=&-2s\func{Im}\int\nolimits_{\Omega }z\partial _{t}(\nabla _{y}\varphi
)\cdot \nabla _{y}\overline{z}dxdydt-2s\func{Im}\int\nolimits_{\Gamma
_{y}}z\partial _{t}\overline{z}\left( \nabla _{y}\varphi \cdot \nu \right)
dS_{y}dxdt  \notag \\
&&+2s\func{Im}\int\nolimits_{\Omega }z\Delta _{y}\varphi \partial _{t}%
\overline{z}dxdydt+2s\func{Im}\int\nolimits_{\Omega }z\partial _{t}(\nabla
_{x}\varphi )\cdot \nabla _{x}\overline{z}dxdydt  \notag \\
&&+2s\func{Im}\int\nolimits_{\Gamma _{x}}z\partial _{t}\overline{z}\left(
\nabla _{x}\varphi \cdot \nu \right) dS_{x}dydt-2s\func{Im}%
\int\nolimits_{\Omega }z\Delta _{x}\varphi \partial _{t}\overline{z}dxdydt. 
\TCItag{3.14}  \label{3.14}
\end{eqnarray}%
In (\ref{3.14}), we used the equality $\func{Re}\left( iz\right) =-\func{Im}%
\left( z\right) $ and $\func{Im}\left( z\right) -\func{Im}\left( \bar{z}%
\right) =2\func{Im}\left( z\right) $, where $\func{Im}\left( z\right) $
denotes the imaginary part of $z\in 
\mathbb{C}
$.%
\begin{eqnarray}
I_{2} &=&-2s\func{Re}\int_{\Omega }i\partial _{t}z\overline{z}\left( \Delta
_{y}\varphi -\Delta _{x}\varphi \right) dxdydt  \notag \\
&=&-2s\func{Im}\int_{\Omega }\partial _{t}\overline{z}z\left( \Delta
_{y}\varphi -\Delta _{x}\varphi \right) dxdydt.  \TCItag{3.15}  \label{3.15}
\end{eqnarray}%
In (\ref{3.15}), we used the equality $\func{Re}\left( iz\right) =\func{Im}%
\left( \bar{z}\right) .$%
\begin{eqnarray}
I_{3} &=&-4s\func{Re}\int_{\Omega }\Delta _{y}z\nabla _{y}\varphi \cdot
\nabla _{y}\overline{z}dxdydt+4s\func{Re}\int_{\Omega }\Delta _{y}z\nabla
_{x}\varphi \cdot \nabla _{x}\overline{z}dxdydt  \notag \\
&=&4s\func{Re}\int\nolimits_{\Omega }\sum\limits_{i,j=1}^{n}\varphi
_{y_{i}y_{j}}\overline{z}_{y_{i}}z_{y_{j}}dxdydt-2s\int\nolimits_{\Omega
}\Delta _{y}\varphi \left\vert \nabla _{y}z\right\vert ^{2}dxdydt  \notag \\
&&+2s\int_{\Gamma _{y}}\left( \partial _{\nu }\varphi \right) \left\vert
\nabla _{y}z\right\vert ^{2}dS_{y}dxdt-4s\func{Re}\int_{\Gamma _{y}}\left(
\partial _{\nu }z\right) \nabla _{y}\varphi \cdot \nabla _{y}\overline{z}%
dS_{y}dxdt  \notag \\
&&-4s\func{Re}\int\nolimits_{\Omega
}\sum\limits_{j=1}^{m}\sum\limits_{i=1}^{n}\varphi _{x_{j}y_{i}}\overline{z}%
_{x_{j}}z_{y_{i}}dxdydt+2s\int\nolimits_{\Omega }\Delta _{x}\varphi
\left\vert \nabla _{y}z\right\vert ^{2}dxdydt  \notag \\
&&-2s\int_{\Gamma _{x}}\left( \partial _{\nu }\varphi \right) \left\vert
\nabla _{y}z\right\vert ^{2}dS_{x}dydt  \notag \\
&&+4s\func{Re}\int_{\Gamma _{y}}\left( \partial _{\nu }z\right) \nabla
_{x}\varphi \cdot \nabla _{x}\overline{z}dS_{y}dxdt.  \TCItag{3.16}
\label{3.16}
\end{eqnarray}%
In (\ref{3.16}), we used the equality $\func{Re}z_{y_{j}}\bar{z}%
_{y_{i}y_{j}}=\frac{1}{2}\left( \left\vert z_{y_{j}}\right\vert ^{2}\right)
_{y_{i}}.$%
\begin{eqnarray}
I_{4} &=&-2s\func{Re}\int_{\Omega }\Delta _{y}z\left( \Delta _{y}\varphi
-\Delta _{x}\varphi \right) \overline{z}dxdydt  \notag \\
&=&-s\func{Re}\int_{\Omega }\Delta _{y}\left( \Delta _{y}\varphi -\Delta
_{x}\varphi \right) \left\vert z\right\vert ^{2}dxdydt  \notag \\
&&+2s\int_{\Omega }\left( \Delta _{y}\varphi -\Delta _{x}\varphi \right)
\left\vert \nabla _{y}z\right\vert ^{2}dxdydt  \notag \\
&&+s\int_{\Gamma _{y}}\partial _{\nu }\left( \Delta _{y}\varphi -\Delta
_{x}\varphi \right) \left\vert z\right\vert ^{2}dS_{y}dxdt  \notag \\
&&-2s\func{Re}\int\nolimits_{\Gamma _{y}}\left( \partial _{\nu }z\right)
\left( \Delta _{y}\varphi -\Delta _{x}\varphi \right) \overline{z}dS_{y}dxdt.
\TCItag{3.17}  \label{3.17}
\end{eqnarray}%
In (\ref{3.17}), we used the equality $\func{Re}\overline{z}\nabla _{y}z=%
\frac{1}{2}\nabla _{y}\left( \left\vert z\right\vert ^{2}\right) .$%
\begin{eqnarray}
I_{5} &=&4s\func{Re}\int_{\Omega }\Delta _{x}z\nabla _{y}\varphi \cdot
\nabla _{y}\overline{z}dxdydt-4s\func{Re}\int_{\Omega }\Delta _{x}z\nabla
_{x}\varphi \cdot \nabla _{x}\overline{z}dxdydt  \notag \\
&=&-4s\func{Re}\int\nolimits_{\Omega
}\sum\limits_{j=1}^{m}\sum\limits_{i=1}^{n}\varphi _{y_{i}x_{j}}\overline{z}%
_{y_{i}}z_{x_{j}}dxdydt-2s\int\nolimits_{\Omega }\sum\limits_{i=1}^{m}\Delta
_{y}\varphi \left\vert \nabla _{x}z\right\vert ^{2}dxdydt  \notag \\
&&-2s\int_{\Gamma _{y}}\left( \partial _{\nu }\varphi \right) \left\vert
\nabla _{x}z\right\vert ^{2}dS_{y}dxdt+4s\func{Re}\int_{\Gamma _{x}}\left(
\partial _{\nu }z\right) \nabla _{y}\varphi \cdot \nabla _{y}\overline{z}%
dS_{x}dydt  \notag \\
&&+4s\func{Re}\int\nolimits_{\Omega }\sum\limits_{i,j=1}^{n}\varphi
_{x_{i}y_{j}}\overline{z}_{x_{i}}z_{x_{j}}dxdydt-2s\int\nolimits_{\Omega
}\Delta _{x}\varphi \left\vert \nabla _{x}z\right\vert ^{2}dxdydt  \notag \\
&&-2s\int_{\Gamma _{x}}\left( \partial _{\nu }\varphi \right) \left\vert
\nabla _{x}z\right\vert ^{2}dS_{x}dydt  \notag \\
&&+4s\func{Re}\int_{\Gamma _{x}}\left( \partial _{\nu }z\right) \nabla
_{x}\varphi \cdot \nabla _{x}\overline{z}dS_{x}dydt.  \TCItag{3.18}
\label{3.18}
\end{eqnarray}%
In (\ref{3.18}), we used the equality $\func{Re}z_{x_{j}}\overline{z}%
_{x_{j}y_{i}}=\frac{1}{2}\left( \left\vert z_{x_{j}}\right\vert ^{2}\right)
_{y_{i}}.$%
\begin{eqnarray}
I_{6} &=&2s\func{Re}\int_{\Omega }\Delta _{x}z\left( \Delta _{y}\varphi
-\Delta _{x}\varphi \right) \overline{z}dxdydt  \notag \\
&=&-2s\int_{\Omega }\left( \Delta _{y}\varphi -\Delta _{x}\varphi \right)
\left\vert \nabla _{x}z\right\vert ^{2}dxdydt  \notag \\
&&+s\int_{\Omega }\Delta _{x}\left( \Delta _{y}\varphi -\Delta _{x}\varphi
\right) )\left\vert z\right\vert ^{2}dxdydt  \notag \\
&&-s\int\nolimits_{\Gamma _{x}}\partial _{\nu }\left( \Delta _{y}\varphi
-\Delta _{x}\varphi \right) \left\vert z\right\vert ^{2}dS_{x}dydt  \notag \\
&&+2s\func{Re}\int\nolimits_{\Gamma _{x}}\left( \partial _{\nu }z\right)
\left( \Delta _{y}\varphi -\Delta _{x}\varphi \right) \overline{z}dS_{x}dydt.
\TCItag{3.19}  \label{3.19}
\end{eqnarray}%
In (\ref{3.19}), we used the equality $\func{Re}\overline{z}$ $\nabla _{x}z=%
\frac{1}{2}\nabla _{x}\left( \left\vert z\right\vert ^{2}\right) .$%
\begin{eqnarray}
I_{7} &=&-4s^{3}\func{Re}\int_{\Omega }\left( \left\vert \nabla _{y}\varphi
\right\vert ^{2}-\left\vert \nabla _{x}\varphi \right\vert ^{2}\right)
z\left( \nabla _{y}\varphi \cdot \nabla _{y}\overline{z}-\nabla _{x}\varphi
\cdot \nabla _{x}\overline{z}\right) dxdydt  \notag \\
&=&2s^{3}\int_{\Omega }\left\vert z\right\vert ^{2}\left( \left\vert \nabla
_{y}\varphi \right\vert ^{2}-\left\vert \nabla _{x}\varphi \right\vert
^{2}\right) \left( \Delta _{y}\varphi -\Delta _{x}\varphi \right) dxdydt 
\notag \\
&&-2s^{3}\int_{\Gamma _{y}}\left( \partial _{\nu }\varphi \right) \left\vert
z\right\vert ^{2}\left( \left\vert \nabla _{y}\varphi \right\vert
^{2}-\left\vert \nabla _{x}\varphi \right\vert ^{2}\right) dS_{y}dxdt  \notag
\\
&&+2s^{3}\int_{\Omega }\left\vert z\right\vert ^{2}\nabla _{y}\varphi \cdot
\nabla _{y}\left( \left\vert \nabla _{y}\varphi \right\vert ^{2}-\left\vert
\nabla _{x}\varphi \right\vert ^{2}\right) dxdydt  \notag \\
&&+2s^{3}\int_{\Gamma _{x}}\partial _{\nu }\varphi \left\vert z\right\vert
^{2}\left( \left\vert \nabla _{y}\varphi \right\vert ^{2}-\left\vert \nabla
_{x}\varphi \right\vert ^{2}\right) dS_{x}dydt  \notag \\
&&-2s^{3}\int_{\Omega }\left\vert z\right\vert ^{2}\nabla _{x}\varphi \cdot
\nabla _{x}\left( \left\vert \nabla _{y}\varphi \right\vert ^{2}-\left\vert
\nabla _{x}\varphi \right\vert ^{2}\right) dxdydt.  \TCItag{3.20}
\label{3.20}
\end{eqnarray}%
In (\ref{3.20}), we used the equality $\func{Re}z\nabla _{y}\overline{z}=%
\frac{1}{2}\left( \nabla _{y}\left\vert z\right\vert ^{2}\right) $ and $%
\func{Re}z\nabla _{x}\overline{z}=\frac{1}{2}\left( \nabla _{x}\left\vert
z\right\vert ^{2}\right) .$%
\begin{eqnarray}
I_{8} &=&-2s^{3}\func{Re}\int_{\Omega }\left( \left\vert \nabla _{y}\varphi
\right\vert ^{2}-\left\vert \nabla _{x}\varphi \right\vert ^{2}\right)
\left( \Delta _{y}\varphi -\Delta _{x}\varphi \right) z\overline{z}dxdydt 
\notag \\
&=&-2s^{3}\int_{\Omega }\left( \left\vert \nabla _{y}\varphi \right\vert
^{2}-\left\vert \nabla _{x}\varphi \right\vert ^{2}\right) \left( \Delta
_{y}\varphi -\Delta _{x}\varphi \right) \left\vert z\right\vert ^{2}dxdydt. 
\TCItag{3.21}  \label{3.21}
\end{eqnarray}%
Hence, we can rewrite (\ref{3.13})%
\begin{equation*}
2\func{Re}\left( P_{s}^{+}z,P_{s}^{-}z\right) _{L^{2}\left( \Omega \right)
}=J_{1}+J_{2}+J_{3}+J_{4}+J_{5}+J_{6}+B_{0},
\end{equation*}%
\bigskip where%
\begin{eqnarray*}
J_{1} &=&-2s\func{Im}\int\nolimits_{\Omega }z\partial _{t}(\nabla
_{y}\varphi )\cdot \nabla _{y}\overline{z}dxdydt+2s\func{Im}%
\int\nolimits_{\Omega }z\partial _{t}(\nabla _{x}\varphi )\cdot \nabla _{x}%
\overline{z}dxdydt, \\
J_{2} &=&4s\func{Re}\sum\limits_{i,j=1}^{n}\int\nolimits_{\Omega }\varphi
_{y_{i}y_{j}}z_{y_{j}}\overline{z}_{y_{i}}dxdydt-4s\func{Re}%
\sum\limits_{j=1}^{m}\sum\limits_{i=1}^{n}\int\nolimits_{\Omega }\varphi
_{y_{i}x_{j}}z_{y_{i}}\overline{z}_{xj}dxdydt, \\
J_{3} &=&-s\int_{\Omega }\left\vert z\right\vert ^{2}\Delta _{y}\left(
\Delta _{y}\varphi -\Delta _{x}\varphi \right) dxdydt, \\
J_{4} &=&-4s\func{Re}\sum\limits_{j=1}^{m}\sum\limits_{i=1}^{n}\int%
\nolimits_{\Omega }\varphi _{y_{i}x_{j}}\overline{z}%
_{y_{i}}z_{x_{j}}dxdydt+4s\func{Re}\sum\limits_{i,j=1}^{n}\int\nolimits_{%
\Omega }\varphi _{x_{i}x_{j}}\overline{z}_{x_{i}}z_{x_{j}}dxdydt,
\end{eqnarray*}%
\begin{eqnarray*}
J_{5} &=&s\int_{\Omega }\left\vert z\right\vert ^{2}\Delta _{x}\left( \Delta
_{y}\varphi -\Delta _{x}\varphi \right) dxdydt, \\
J_{6} &=&2s^{3}\int_{\Omega }\left\vert z\right\vert ^{2}\nabla _{y}\varphi
\cdot \nabla _{y}\left( \left\vert \nabla _{y}\varphi \right\vert
^{2}-\left\vert \nabla _{x}\varphi \right\vert ^{2}\right) dxdydt \\
&&-2s^{3}\int_{\Omega }\left\vert z\right\vert ^{2}\nabla _{x}\varphi \cdot
\nabla _{x}\left( \left\vert \nabla _{y}\varphi \right\vert ^{2}-\left\vert
\nabla _{x}\varphi \right\vert ^{2}\right) dxdydt
\end{eqnarray*}%
and%
\begin{eqnarray*}
B_{0} &=&-2s\func{Im}\int\nolimits_{\Gamma _{y}}z\partial _{t}\overline{z}%
\left( \nabla _{y}\varphi \cdot \nu \right) dS_{y}dxdt \\
&&+4s\func{Re}\int\nolimits_{\Gamma _{y}}(\partial _{\nu }z)\left( \nabla
_{x}\varphi \cdot \nabla _{x}\overline{z}-\nabla _{y}\varphi \cdot \nabla
_{y}\overline{z}\right) dS_{y}dxdt \\
&&+2s\int\nolimits_{\Gamma _{y}}(\partial _{\nu }\varphi )\left( \left\vert
\nabla _{y}z\right\vert ^{2}-\left\vert \nabla _{x}z\right\vert ^{2}\right)
dS_{y}dxdt \\
&&-2s^{3}\int_{\Gamma _{y}}(\partial _{\nu }\varphi )\left\vert z\right\vert
^{2}\left( \left\vert \nabla _{y}\varphi \right\vert ^{2}-\left\vert \nabla
_{x}\varphi \right\vert ^{2}\right) dS_{y}dxdt \\
&&-2s\int\nolimits_{\Gamma _{y}}(\partial _{\nu }z)\left( \Delta _{y}\varphi
-\Delta _{x}\varphi \right) \overline{z}dS_{y}dxdt+s\int\nolimits_{\Gamma
_{y}}\partial _{\nu }\left( \Delta _{y}\varphi -\Delta _{x}\varphi \right)
\left\vert z\right\vert ^{2}dS_{y}dxdt \\
&&+2s\func{Im}\int\nolimits_{\Gamma _{x}}z\partial _{t}\overline{z}\left(
\nabla _{x}\varphi \cdot \nu \right) dS_{x}dydt \\
&&+4s\func{Re}\int_{\Gamma _{x}}(\partial _{\nu }z)\left( \nabla _{y}\varphi
\cdot \nabla _{y}\overline{z}-\nabla _{x}\varphi \cdot \nabla _{x}\overline{z%
}\right) dS_{x}dydt \\
&&-2s\int\nolimits_{\Gamma _{x}}(\partial _{\nu }\varphi )\left( \left\vert
\nabla _{y}z\right\vert ^{2}-\left\vert \nabla _{x}z\right\vert ^{2}\right)
dS_{x}dydt \\
&&+2s^{3}\int_{\Gamma _{x}}(\partial _{\nu }\varphi )\left\vert z\right\vert
^{2}\left( \left\vert \nabla _{y}\varphi \right\vert ^{2}-\left\vert \nabla
_{x}\varphi \right\vert ^{2}\right) dS_{x}dydt \\
&&+2s\int_{\Gamma _{x}}(\partial _{\nu }z)\left( \Delta _{y}\varphi -\Delta
_{x}\varphi \right) \overline{z}dS_{x}dydt-s\int_{\Gamma _{x}}\partial _{\nu
}\left( \Delta _{y}\varphi -\Delta _{x}\varphi \right) \left\vert
z\right\vert ^{2}dS_{x}dydt.
\end{eqnarray*}%
Next, we shall estimate $J_{k},$ $1\leq k\leq 6$ and $B_{0}$ using the
following elemantary properties of the weight function:%
\begin{equation*}
\begin{array}{cc}
\partial _{t}\varphi =\left( -2\gamma \beta t\right) \varphi , & \varphi
_{x_{i}x_{i}}=\gamma \varphi \left( 2+\gamma \psi _{x_{i}}^{2}\right) , \\ 
\varphi _{x_{i}y_{j}}=\gamma ^{2}\varphi \psi _{x_{i}}\psi _{y_{j}}, & 
\varphi _{x_{i}x_{j}}=\gamma \varphi \left( \psi _{x_{i}x_{j}}+\gamma \psi
_{x_{i}}\psi _{x_{j}}\right) , \\ 
\varphi _{y_{i}y_{j}}=\gamma \varphi \left( \psi _{y_{i}y_{j}}+\gamma \psi
_{y_{i}}\psi _{y_{j}}\right) , & \nabla _{x}\varphi =\gamma \varphi \nabla
_{x}\psi , \\ 
\nabla _{y}\varphi =\gamma \varphi \nabla _{y}\psi , & \partial _{t}(\nabla
_{x}\varphi )=\left( -2\gamma ^{2}\beta t\right) \varphi \nabla _{x}\psi ,
\\ 
\partial _{t}(\nabla _{y}\varphi )=\left( -2\gamma ^{2}\beta t\right)
\varphi \nabla _{y}\psi , & \Delta _{x}\varphi =\gamma \varphi \left( \Delta
_{x}\psi +\gamma \left\vert \nabla _{x}\psi \right\vert ^{2}\right) , \\ 
\Delta _{y}\varphi =\gamma \varphi \left( \Delta _{y}\psi +\gamma \left\vert
\nabla _{y}\psi \right\vert ^{2}\right) , & \Delta _{y}\varphi -\Delta
_{x}\varphi =\gamma \varphi d_{1}\left( \psi \right) +\gamma ^{2}\varphi
d_{2}\left( \psi \right) ,%
\end{array}%
\end{equation*}%
where%
\begin{eqnarray*}
d_{1}\left( \psi \right) &=&\Delta _{y}\psi -\Delta _{x}\psi , \\
d_{2}\left( \psi \right) &=&\left\vert \nabla _{y}\psi \right\vert
^{2}-\left\vert \nabla _{x}\psi \right\vert ^{2}.
\end{eqnarray*}%
Then, we obtain%
\begin{eqnarray}
J_{1} &=&-2s\func{Im}\int\nolimits_{\Omega }z\partial _{t}(\nabla
_{y}\varphi )\cdot \nabla _{y}\overline{z}dxdydt+2s\func{Im}%
\int\nolimits_{\Omega }z\partial _{t}(\nabla _{x}\varphi )\cdot \nabla _{x}%
\overline{z}dxdydt  \notag \\
&=&-2s\func{Im}\int\nolimits_{\Omega }\left( -2\gamma ^{2}\beta t\right)
z\varphi \nabla _{y}\psi \cdot \nabla _{y}\overline{z}dxdydt  \notag \\
&&+2s\func{Im}\int\nolimits_{\Omega }\left( -2\gamma ^{2}\beta t\right)
z\varphi \nabla _{x}\psi \cdot \nabla _{x}\overline{z}dxdydt  \TCItag{3.22}
\label{3.22}
\end{eqnarray}%
and%
\begin{eqnarray}
J_{2} &=&\sum\limits_{i,j=1}^{m}\func{Re}\int\nolimits_{\Omega }4s\varphi
_{y_{i}y_{j}}z_{y_{j}}\overline{z}_{y_{i}}dxdydt-\func{Re}%
\sum\limits_{j=1}^{m}\sum\limits_{i=1}^{n}\int\nolimits_{\Omega }4sz_{y_{i}}%
\overline{z}_{x_{j}}\varphi _{y_{i}x_{j}}dxdydt  \notag \\
&=&\sum\limits_{i,j=1}^{m}\func{Re}\int\nolimits_{\Omega }4s\gamma \varphi
\left( \psi _{y_{i}y_{j}}+\gamma \psi _{y_{i}}\psi _{y_{j}}\right) z_{y_{j}}%
\overline{z}_{y_{i}}dxdydt+\int\nolimits_{\Omega }4s\gamma ^{2}\varphi
\left\vert \nabla _{y}\psi \cdot \nabla _{y}z\right\vert ^{2}dxdydt  \notag
\\
&&-\sum\limits_{j=1}^{m}\sum\limits_{i=1}^{n}\int\nolimits_{\Omega }4s\gamma
^{2}\varphi \left( \nabla _{y}\psi \cdot \nabla _{y}z\right) \left( \nabla
_{x}\psi \cdot \nabla _{x}\overline{z}\right) dxdydt.  \TCItag{3.23}
\label{3.23}
\end{eqnarray}%
Before estimating $J_{3}$, we can directly verify%
\begin{eqnarray*}
\Delta _{y}\left( \varphi d_{2}\left( \psi \right) \right) &=&\left( \Delta
_{y}\varphi \right) d_{2}\left( \psi \right) +2\nabla _{y}\varphi \cdot
\nabla _{y}\left( d_{2}\left( \psi \right) \right) +\varphi \Delta
_{y}\left( d_{2}\left( \psi \right) \right) \\
&=&\gamma \varphi \left( \Delta _{y}\psi \right) d_{2}\left( \psi \right)
+\gamma ^{2}\varphi \left\vert \nabla _{y}\psi \right\vert ^{2}d_{2}\psi
+2\gamma \varphi \nabla _{y}\psi \cdot \nabla _{y}\left( d_{2}\left( \psi
\right) \right) \\
&&+\varphi \Delta _{y}\left( d_{2}\left( \psi \right) \right) , \\
\Delta _{y}\left( \varphi d_{1}\left( \psi \right) \right) &=&\left( \Delta
_{y}\varphi \right) d_{1}\left( \psi \right) +2\nabla _{y}\varphi \cdot
\nabla _{y}\left( d_{1}\left( \psi \right) \right) +\varphi \Delta
_{y}\left( d_{1}\left( \psi \right) \right) \\
&=&\gamma \varphi \left( \Delta _{y}\psi \right) d_{1}\left( \psi \right)
+\gamma ^{2}\varphi \left\vert \nabla _{y}\psi \right\vert ^{2}d_{1}\psi
+2\gamma \varphi \nabla _{y}\psi \cdot \nabla _{y}\left( d_{1}\left( \psi
\right) \right) \\
&&+\varphi \Delta _{y}\left( d_{1}\left( \psi \right) \right) .
\end{eqnarray*}%
Then, we have%
\begin{eqnarray}
J_{3} &=&-s\int_{\Omega }\left\vert z\right\vert ^{2}\Delta _{y}\left(
\Delta _{y}\varphi -\Delta _{x}\varphi \right) dxdydt  \notag \\
&=&-\int_{\Omega }s\gamma ^{2}\varphi \left\vert z\right\vert ^{2}\left(
d_{1}\left( \psi \right) \left( \Delta _{y}\psi \right) +\Delta _{y}\left(
d_{2}\left( \psi \right) \right) \right) dxdydt  \notag \\
&&-\int_{\Omega }s\gamma ^{3}\varphi \left\vert z\right\vert ^{2}\left(
d_{1}\left( \psi \right) \left\vert \nabla _{y}\psi \right\vert ^{2}+\left(
\Delta _{y}\psi \right) d_{2}\left( \psi \right) +2\nabla _{y}\psi \cdot
\nabla _{y}\left( d_{2}\left( \psi \right) \right) \right) dxdydt  \notag \\
&&-\int_{\Omega }s\gamma ^{4}\varphi \left\vert z\right\vert ^{2}\left\vert
\nabla _{y}\psi \right\vert ^{2}d_{2}\left( \psi \right) dxdydt, 
\TCItag{3.24}  \label{3.24}
\end{eqnarray}%
\begin{eqnarray}
J_{4} &=&-\sum\limits_{j=1}^{m}\sum\limits_{i=1}^{n}\func{Re}%
\int\nolimits_{\Omega }4s\varphi _{y_{i}x_{j}}\overline{z}%
_{y_{i}}z_{x_{j}}dxdydt+\sum\limits_{i,j=1}^{n}\func{Re}\int\nolimits_{%
\Omega }4s\varphi _{x_{i}x_{j}}\overline{z}_{x_{i}}z_{x_{j}}dxdydt  \notag \\
&=&-\func{Re}\int\nolimits_{\Omega }4s\gamma ^{2}\varphi \left( \nabla
_{y}\psi \cdot \nabla _{y}\overline{z}\right) \left( \nabla _{x}\psi \cdot
\nabla _{x}z\right) dxdydt  \notag \\
&&+\sum\limits_{i,j=1}^{n}\func{Re}\int\nolimits_{\Omega }4s\gamma \varphi
\psi _{x_{i}x_{j}}\overline{z}_{x_{i}}z_{x_{j}}dxdydt  \notag \\
&&+\int\nolimits_{\Omega }4s\gamma ^{2}\varphi \left\vert \nabla _{x}\psi
\cdot \nabla _{x}z\right\vert ^{2}dxdydt.  \TCItag{3.25}  \label{3.25}
\end{eqnarray}%
Since%
\begin{eqnarray*}
\Delta _{x}\left( \varphi d_{2}\left( \psi \right) \right) &=&\left( \Delta
_{x}\varphi \right) d_{2}\left( \psi \right) +2\nabla _{x}\varphi \cdot
\nabla _{x}\left( d_{2}\left( \psi \right) \right) +\varphi \Delta
_{x}\left( d_{2}\left( \psi \right) \right) \\
&=&\gamma \varphi \left( \Delta _{x}\psi \right) d_{2}\left( \psi \right)
+\gamma ^{2}\varphi \left\vert \nabla _{x}\psi \right\vert ^{2}d_{2}\psi
+2\gamma \varphi \nabla _{x}\psi \cdot \nabla _{x}\left( d_{2}\left( \psi
\right) \right) \\
&&+\varphi \Delta _{x}\left( d_{2}\left( \psi \right) \right) ,
\end{eqnarray*}%
we see that%
\begin{eqnarray}
J_{5} &=&s\int_{\Omega }\left\vert z\right\vert ^{2}\Delta _{x}\left( \Delta
_{y}\varphi -\Delta _{x}\varphi \right) dxdydt  \notag \\
&=&\int_{\Omega }s\gamma ^{2}\varphi \left\vert z\right\vert ^{2}\left(
d_{1}\left( \psi \right) \left( \Delta _{x}\psi \right) +\Delta _{x}\left(
d_{2}\left( \psi \right) \right) \right) dxdydt  \notag \\
&&+\int_{\Omega }s\gamma ^{3}\varphi \left\vert z\right\vert ^{2}\left(
d_{1}\left( \psi \right) \left\vert \nabla _{x}\psi \right\vert ^{2}+\left(
\Delta _{x}\psi \right) d_{2}\left( \psi \right) +2\nabla _{x}\psi \cdot
\nabla _{x}\left( d_{2}\left( \psi \right) \right) \right) dxdydt  \notag \\
&&+\int_{\Omega }s\gamma ^{4}\varphi \left\vert z\right\vert ^{2}\left\vert
\nabla _{x}\psi \right\vert ^{2}d_{2}\left( \psi \right) dxdydt. 
\TCItag{3.26}  \label{3.26}
\end{eqnarray}%
Since%
\begin{equation*}
\nabla _{y}\varphi \cdot \nabla _{y}\left( \left\vert \nabla _{y}\varphi
\right\vert ^{2}-\left\vert \nabla _{x}\varphi \right\vert ^{2}\right)
=2\gamma ^{4}\varphi ^{3}d_{2}\left( \psi \right) \nabla _{y}\psi \cdot
\nabla _{y}\psi +\gamma ^{3}\varphi ^{3}\nabla _{y}\psi \cdot \nabla
_{y}\left( d_{2}\left( \psi \right) \right)
\end{equation*}%
and%
\begin{equation*}
\nabla _{x}\varphi \cdot \nabla _{x}\left( \left\vert \nabla _{y}\varphi
\right\vert ^{2}-\left\vert \nabla _{x}\varphi \right\vert ^{2}\right)
=2\gamma ^{4}\varphi ^{3}d_{2}\left( \psi \right) \nabla _{x}\psi \cdot
\nabla _{x}\psi +\gamma ^{3}\varphi ^{3}\nabla _{x}\psi \cdot \nabla
_{x}\left( d_{2}\left( \psi \right) \right) ,
\end{equation*}%
we have 
\begin{eqnarray}
J_{6} &=&2s^{3}\int_{\Omega }\left\vert z\right\vert ^{2}\nabla _{y}\varphi
\cdot \nabla _{y}\left( \left\vert \nabla _{y}\varphi \right\vert
^{2}-\left\vert \nabla _{x}\varphi \right\vert ^{2}\right) dxdydt  \notag \\
&&-2s^{3}\int_{\Omega }\left\vert z\right\vert ^{2}\nabla _{x}\varphi \cdot
\nabla _{x}\left( \left\vert \nabla _{y}\varphi \right\vert ^{2}-\left\vert
\nabla _{x}\varphi \right\vert ^{2}\right) dxdydt  \notag \\
&=&4\int_{\Omega }s^{3}\gamma ^{4}\varphi ^{3}\left\vert z\right\vert
^{2}\left( d_{2}\left( \psi \right) \right) ^{2}dxdydt  \notag \\
&&+\int_{\Omega }2s^{3}\gamma ^{3}\varphi ^{3}\left\vert z\right\vert
^{2}\nabla _{y}\psi \cdot \nabla _{y}\left( d_{2}\left( \psi \right) \right)
dxdydt  \notag \\
&&-\int_{\Omega }2s^{3}\gamma ^{3}\varphi ^{3}\left\vert z\right\vert
^{2}\nabla _{x}\psi \cdot \nabla _{x}\left( d_{2}\left( \psi \right) \right)
dxdydt.  \TCItag{3.27}  \label{3.27}
\end{eqnarray}%
Finally, the boundary term is obtained as follows:%
\begin{eqnarray*}
B_{0} &=&-2\func{Im}\int\nolimits_{\Gamma _{y}}s\gamma \varphi z\partial _{t}%
\overline{z}\left( \nabla _{y}\psi \cdot \nu \right) dS_{y}dxdt \\
&&+4\func{Re}\int\nolimits_{\Gamma _{y}}s\gamma \varphi \left( \partial
_{\nu }z\right) \left( \nabla _{x}\psi \cdot \nabla _{x}\overline{z}-\nabla
_{y}\psi \cdot \nabla _{y}\overline{z}\right) dS_{y}dxdt \\
&&+\int\nolimits_{\Gamma _{y}}2s\gamma \varphi \left( \partial _{\nu }\psi
\right) \left( \left\vert \nabla _{y}z\right\vert ^{2}-\left\vert \nabla
_{x}z\right\vert ^{2}\right) dS_{y}dxdt \\
&&-\int_{\Gamma _{y}}2s^{3}\gamma ^{3}\varphi ^{3}\partial _{v}\psi
\left\vert z\right\vert ^{2}d_{2}\left( \psi \right) dS_{y}dxdt \\
&&-\int\nolimits_{\Gamma _{y}}2s\left( \gamma \varphi d_{1}\left( \psi
\right) +\gamma ^{2}\varphi d_{2}\left( \psi \right) \right) \overline{z}%
\left( \partial _{\nu }z\right) dS_{y}dxdt \\
&&+\int\nolimits_{\Gamma _{y}}s\left( \left( \gamma ^{2}\varphi d_{1}\left(
\psi \right) +\gamma ^{3}\varphi d_{2}\left( \psi \right) \right) \left(
\partial _{\nu }\psi \right) \right. \\
&&\left. +\gamma ^{2}\varphi \left( \partial _{\nu }\left( d_{2}\left( \psi
\right) \right) \right) \right) \left\vert z\right\vert ^{2}dS_{y}dxdt \\
&&+2\func{Im}\int\nolimits_{\Gamma _{x}}s\gamma \varphi z\partial _{t}%
\overline{z}\left( \nabla _{x}\psi \cdot \nu \right) dS_{x}dydt \\
&&+4\func{Re}\int_{\Gamma _{x}}s\gamma \varphi \left( \partial _{v}z\right)
\left( \nabla _{y}\psi \cdot \nabla _{y}\overline{z}-\nabla _{x}\psi \cdot
\nabla _{x}\overline{z}\right) dS_{x}dydt
\end{eqnarray*}%
\begin{eqnarray}
&&-\int\nolimits_{\Gamma _{x}}2s\gamma \varphi \left( \partial _{\nu }\psi
\right) \left( \left\vert \nabla _{y}z\right\vert ^{2}-\left\vert \nabla
_{x}z\right\vert ^{2}\right) dS_{x}dydt  \notag \\
&&+\int_{\Gamma _{x}}2s^{3}\gamma ^{3}\varphi ^{3}\left( \partial _{\nu
}\psi \right) \left\vert z\right\vert ^{2}d_{2}\left( \psi \right) dS_{x}dydt
\notag \\
&&+\int_{\Gamma _{x}}2s\left( \gamma \varphi d_{1}\left( \psi \right)
+\gamma ^{2}\varphi d_{2}\left( \psi \right) \right) \overline{z}\left(
\partial _{\nu }z\right) dS_{x}dydt  \notag \\
&&-\int_{\Gamma _{x}}s\left( \left( \gamma ^{2}\varphi d_{1}\left( \psi
\right) +\gamma ^{3}\varphi d_{2}\left( \psi \right) \right) \left( \partial
_{\nu }\psi \right) \right.  \notag \\
&&\left. +\gamma ^{2}\varphi \partial _{\nu }\left( d_{2}\left( \psi \right)
\right) \right) \left\vert z\right\vert ^{2}dS_{x}dydt.  \TCItag{3.28}
\label{3.28}
\end{eqnarray}%
Then from (\ref{3.22})-(\ref{3.28}) we can write%
\begin{eqnarray*}
2\func{Re}\left( P_{s}^{+}z,P_{s}^{-}z\right) _{L^{2}\left( \Omega \right) }
&=&\sum\limits_{i,j=1}^{m}\func{Re}\int\nolimits_{\Omega }4s\gamma \varphi
\psi _{y_{i}y_{j}}z_{y_{j}}\overline{z}_{y_{i}}dxdydt \\
&&+\sum\limits_{i,j=1}^{n}\func{Re}\int\nolimits_{\Omega }4s\gamma \varphi
\psi _{x_{i}x_{j}}\overline{z}_{x_{i}}z_{x_{j}}dxdydt \\
&&+\int\nolimits_{\Omega }4s\gamma ^{2}\varphi \left\vert \nabla _{y}\psi
\cdot \nabla _{y}z-\nabla _{x}\psi \cdot \nabla _{x}z\right\vert ^{2}dxdydt
\\
&&+\int\nolimits_{\Omega }4s^{3}\gamma ^{4}\varphi ^{3}\left\vert
z\right\vert ^{2}\left( d_{2}\left( \psi \right) \right)
^{2}dxdydt+B_{0}+X_{1}+X_{2},
\end{eqnarray*}%
where%
\begin{equation*}
X_{1}=\int\nolimits_{\Omega }2s^{3}\gamma ^{3}\varphi ^{3}\left\vert
z\right\vert ^{2}d_{5}\left( \psi \right) dxdydt,
\end{equation*}%
\begin{eqnarray*}
X_{2} &=&-2\func{Im}\int\nolimits_{\Omega }s\left( -2\gamma ^{2}\beta
t\right) z\varphi \nabla _{y}\psi \cdot \nabla _{y}\overline{z}dxdydt \\
&&+2\func{Im}\int\nolimits_{\Omega }s\left( -2\gamma ^{2}\beta t\right)
z\varphi \nabla _{x}\psi \cdot \nabla _{x}\overline{z}dxdydt \\
&&-\int\nolimits_{\Omega }s\gamma ^{4}\varphi \left\vert z\right\vert
^{2}\left( d_{2}\left( \psi \right) \right) ^{2}dxdydt \\
&&-\int\nolimits_{\Omega }s\gamma ^{2}\varphi \left\vert z\right\vert
^{2}d_{3}\left( \psi \right) dxdydt \\
&&-\int\nolimits_{\Omega }s\gamma ^{3}\varphi \left\vert z\right\vert
^{2}d_{4}\left( \psi \right) dxdydt,
\end{eqnarray*}%
and%
\begin{eqnarray*}
d_{3} &:&=d_{3}\left( \psi \right) =\left( d_{1}\left( \psi \right) \right)
^{2}+\Delta _{y}\left( d_{2}\left( \psi \right) \right) -\Delta _{x}\left(
d_{2}\left( \psi \right) \right) , \\
d_{4} &:&=d_{4}\left( \psi \right) =d_{1}\left( \psi \right) d_{2}\left(
\psi \right) +\nabla _{y}\psi \cdot \nabla _{y}\left( d_{2}\left( \psi
\right) \right) -\nabla _{x}\psi \cdot \nabla _{x}\left( d_{2}\left( \psi
\right) \right) , \\
d_{5} &:&=d_{5}\left( \psi \right) =\nabla _{y}\psi \cdot \nabla _{y}\left(
d_{2}\left( \psi \right) \right) -\nabla _{x}\psi \cdot \nabla _{x}\left(
d_{2}\left( \psi \right) \right) .
\end{eqnarray*}%
Since 
\begin{equation*}
\int\nolimits_{\Omega }4s\gamma ^{2}\varphi \left\vert \nabla _{y}\psi \cdot
\nabla _{y}z-\nabla _{x}\psi \cdot \nabla _{x}z\right\vert ^{2}dxdydt\geq 0
\end{equation*}%
and for $0<\alpha ,\beta <1,$%
\begin{equation*}
\left\vert x-x_{0}\right\vert ^{2}-\alpha ^{2}|y|^{2}-\beta ^{2}\left\vert
t\right\vert ^{2}\geq \delta _{0}^{2},
\end{equation*}%
we have%
\begin{equation*}
d_{2}^{2}=16(|x-x_{0}|^{2}-\alpha ^{2}|y-y_{0}|^{2})^{2}\geq
16(|x-x_{0}|^{2}-\alpha ^{2}|y-y_{0}|^{2}-\beta ^{2}\left\vert t\right\vert
^{2})^{2}\geq 16\delta _{0}^{2}.
\end{equation*}%
Then, we see that%
\begin{eqnarray}
2\func{Re}\left( P_{s}^{+}z,P_{s}^{-}z\right) _{L^{2}\left( \Omega \right) }
&\geq &-\int\nolimits_{\Omega }8s\alpha \gamma \varphi |\nabla
_{y}z|^{2}dxdydt+\int_{\Omega }8s\gamma \varphi |\nabla _{x}z|^{2}dxdydt 
\notag \\
&&+64\delta _{0}^{2}\int_{\Omega }s^{3}\gamma ^{4}\varphi
^{3}|z|^{2}dxdydt+B_{0}+X_{1}+X_{2}.  \TCItag{3.29}  \label{3.29}
\end{eqnarray}%
Since the signs of the terms of $|\nabla _{x}z|^{2}$ and $|\nabla _{y}z|^{2}$
are different, we need to perform another estimation:%
\begin{eqnarray*}
2\func{Re}\left( P_{s}^{+}z+P_{s}^{-}z,\varphi z\right) _{L^{2}\left( \Omega
\right) } &=&2\func{Re}\int\nolimits_{\Omega }i\partial _{t}z\overline{z}%
\varphi dxdydt+2\func{Re}\int\nolimits_{\Omega }\Delta _{y}z\overline{z}%
\varphi dxdydt \\
&&-2\func{Re}\int\nolimits_{\Omega }\Delta _{x}z\overline{z}\varphi dxdydt \\
&&+2\func{Re}\int\nolimits_{\Omega }s^{2}(\left\vert \nabla _{y}\varphi
\right\vert ^{2}-\left\vert \nabla _{x}\varphi \right\vert ^{2})\varphi z%
\overline{z}dxdydt \\
&&-4\func{Re}\int\nolimits_{\Omega }s\left( \nabla _{y}\varphi \cdot \nabla
_{y}z-\nabla _{x}\varphi \cdot \nabla _{x}z\right) \varphi \overline{z}dxdydt
\\
&&-2\func{Re}\int\nolimits_{\Omega }s\left( \Delta _{y}\varphi -\Delta
_{x}\varphi \right) z\varphi \overline{z}dxdydt \\
&=&K_{1}+K_{2}+K_{3}+K_{4}+K_{5}+K_{6}.
\end{eqnarray*}%
Now we calculate the terms $K_{j},$ $1\leq j\leq 6$ as follows:%
\begin{eqnarray*}
K_{1} &=&2\func{Re}\int\nolimits_{\Omega }i\partial _{t}z\overline{z}\varphi
dxdydt \\
&=&-2\func{Im}\underset{\Omega }{\int }\partial _{t}z\overline{z}\varphi
dxdydt,
\end{eqnarray*}%
\begin{eqnarray*}
K_{2} &=&2\func{Re}\int\nolimits_{\Omega }\Delta _{y}z\varphi \overline{z}%
dxdydt \\
&=&\int\nolimits_{\Omega }\gamma \varphi \left\vert z\right\vert ^{2}\Delta
_{y}\psi dxdydt+\int\nolimits_{\Omega }\gamma ^{2}\varphi \left\vert
z\right\vert ^{2}\left\vert \nabla _{y}\psi \right\vert ^{2}dxdydt \\
&&-2\int\nolimits_{\Gamma _{y}}\gamma \varphi \left( \partial _{\nu }\psi
\right) \left\vert z\right\vert ^{2}dS_{y}dxdt-2\int\nolimits_{\Omega
}\varphi \left\vert \nabla _{y}z\right\vert ^{2}dxdydt \\
&&+2\func{Re}\int\nolimits_{\Gamma _{y}}\left( \partial _{\nu }z\right)
\left( \varphi \overline{z}\right) dS_{y}dxdt,
\end{eqnarray*}%
\begin{eqnarray*}
K_{3} &=&-2\func{Re}\int\nolimits_{\Omega }\Delta _{x}z\overline{z}\varphi
dxdydt \\
&=&2\int\nolimits_{\Omega }\varphi \left\vert \nabla _{x}z\right\vert
^{2}dxdydt-\int\nolimits_{\Omega }\gamma \varphi \Delta _{x}\psi \left\vert
z\right\vert ^{2}dxdydt \\
&&-\int\nolimits_{\Omega }\gamma ^{2}\varphi \left\vert \nabla _{x}\psi
\right\vert ^{2}\left\vert z\right\vert ^{2}dxdydt+\int_{\Gamma _{x}}\gamma
\varphi \left( \partial _{\nu }\psi \right) \left\vert z\right\vert
^{2}dS_{x}dydt \\
&&-2\func{Re}\int_{\Gamma _{x}}\left( \partial _{\nu }z\right) \left(
\varphi \overline{z}\right) dS_{x}dydt,
\end{eqnarray*}%
\begin{eqnarray*}
K_{4} &=&2\func{Re}\int\nolimits_{\Omega }s^{2}\left( \left\vert \nabla
_{y}\varphi \right\vert ^{2}-\left\vert \nabla _{x}\varphi \right\vert
^{2}\right) \varphi z\overline{z}dxdydt \\
&=&2\int\nolimits_{\Omega }s^{2}\gamma ^{2}\varphi ^{3}\left( \left\vert
\nabla _{y}\psi \right\vert ^{2}-\left\vert \nabla _{x}\psi \right\vert
^{2}\right) \left\vert z\right\vert ^{2}dxdydt \\
&=&2\int\nolimits_{\Omega }s^{2}\gamma ^{2}\varphi ^{3}d_{2}\left( \psi
\right) \left\vert z\right\vert ^{2}dxdydt,
\end{eqnarray*}%
\begin{eqnarray*}
K_{5} &=&-4\func{Re}\text{ }\int\nolimits_{\Omega }s\varphi \overline{z}%
\left( \nabla _{y}\varphi \cdot \nabla _{y}z-\nabla _{x}\varphi \cdot \nabla
_{x}z\right) \varphi dxdydt \\
&=&2\int\nolimits_{\Omega }s\varphi ^{2}\gamma \left\vert z\right\vert
^{2}d_{1}\left( \psi \right) dxdydt+4\int\nolimits_{\Omega }s\gamma
^{2}\varphi ^{2}\left\vert z\right\vert ^{2}d_{2}\left( \psi \right) dxdydt
\\
&&-2\int_{\Gamma _{y}}s\gamma \varphi ^{2}\left( \partial _{\nu }\psi
\right) \left\vert z\right\vert ^{2}dS_{y}dxdt+2\int_{\Gamma _{x}}s\gamma
\varphi ^{2}\left( \partial _{\nu }\psi \right) \left\vert z\right\vert
^{2}dS_{x}dydt,
\end{eqnarray*}%
\begin{eqnarray*}
K_{6} &=&-2\func{Re}\int\nolimits_{\Omega }s\varphi z\overline{z}\left(
\Delta _{y}\varphi -\Delta _{x}\varphi \right) dxdydt \\
&=&-2\int\nolimits_{\Omega }s\gamma \varphi ^{2}\left( \Delta _{y}\psi
-\Delta _{x}\psi \right) \left\vert z\right\vert ^{2}dxdydt \\
&&-2\int\nolimits_{\Omega }s\gamma ^{2}\varphi ^{2}\left\vert z\right\vert
^{2}\left( \left\vert \nabla _{y}\psi \right\vert ^{2}-\left\vert \nabla
_{x}\psi \right\vert ^{2}\right) dxdydt \\
&=&-2\int\nolimits_{\Omega }s\gamma \varphi ^{2}d_{1}\left( \psi \right)
\left\vert z\right\vert ^{2}dxdydt-2\int\nolimits_{\Omega }s\gamma
^{2}\varphi ^{2}\left\vert z\right\vert ^{2}d_{2}\left( \psi \right) dxdydt.
\end{eqnarray*}%
Then we obtain%
\begin{eqnarray}
2\func{Re}\int\nolimits_{\Omega }\left( P_{s}^{+}z+P_{s}^{-}z\right) \varphi 
\overline{z}dxdydt &=&-2\underset{\Omega }{\int }\varphi \left\vert \nabla
_{y}z\right\vert ^{2}dxdydt+2\underset{\Omega }{\int }\varphi \left\vert
\nabla _{x}z\right\vert ^{2}dxdydt  \notag \\
&&+B_{1}+X_{3}+X_{4},  \TCItag{3.30}  \label{3.30}
\end{eqnarray}%
where%
\begin{eqnarray*}
X_{3} &=&2\int\nolimits_{\Omega }s^{2}\gamma ^{2}\varphi ^{3}d_{2}\left(
\psi \right) \left\vert z\right\vert ^{2}dxdydt, \\
X_{4} &=&-2\func{Im}\underset{\Omega }{\int }\partial _{t}z\overline{z}%
\varphi dxdydt+\int\nolimits_{\Omega }\gamma \varphi \left\vert z\right\vert
^{2}d_{1}\left( \psi \right) dxdydt \\
&&+\int\nolimits_{\Omega }s\gamma ^{2}\varphi \left\vert z\right\vert
^{2}d_{2}\left( \psi \right) dxdydt+2\int\nolimits_{\Omega }s\gamma
^{2}\varphi ^{2}\left\vert z\right\vert ^{2}d_{2}\left( \psi \right) dxdydt.
\end{eqnarray*}%
We note that%
\begin{eqnarray*}
B_{1} &=&-\int_{\Gamma _{y}}\gamma \varphi \left( \partial _{\nu }\psi
\right) \left\vert z\right\vert ^{2}dS_{y}dxdt+2\func{Re}\int_{\Gamma
_{y}}\left( \partial _{\nu }z\right) \left( \varphi \overline{z}\right)
dS_{y}dxdt \\
&&-2\int_{\Gamma _{y}}s\gamma \varphi ^{2}\left( \partial _{\nu }\psi
\right) \left\vert z\right\vert ^{2}dS_{y}dxdt+\int_{\Gamma _{x}}\gamma
\varphi \left( \partial _{\nu }\psi \right) \left\vert z\right\vert
^{2}dS_{x}dydt \\
&&-2\func{Re}\int_{\Gamma _{x}}\left( \partial _{\nu }z\right) \left(
\varphi \overline{z}\right) dS_{x}dydt+2\int_{\Gamma _{x}}s\gamma \varphi
^{2}\left( \partial _{\nu }\psi \right) \left\vert z\right\vert
^{2}dS_{x}dxdt=0,
\end{eqnarray*}%
since $z=0$ on $\Gamma _{x}\cup \Gamma _{y}$. We multiply (\ref{3.30}) by $%
-s\gamma \left( 4\alpha +\mu \right) $, then we have%
\begin{eqnarray}
-2\func{Re}\int\nolimits_{\Omega }(4\alpha +\mu )\left(
P_{s}^{+}z+P_{s}^{-}z\right) s\gamma \varphi \overline{z}dxdydt &=&2\underset%
{\Omega }{\int }\left( 4\alpha +\mu \right) s\gamma \varphi \left\vert
\nabla _{y}z\right\vert ^{2}dxdydt  \notag \\
&&-2\underset{\Omega }{\int }\left( 4\alpha +\mu \right) s\gamma \varphi
\left\vert \nabla _{x}z\right\vert ^{2}dxdydt  \notag \\
&&+X_{5}+X_{6},  \TCItag{3.31}  \label{3.31}
\end{eqnarray}%
where we choose $%
\mu
>0$ later and 
\begin{eqnarray*}
X_{5} &=&-2\int\nolimits_{\Omega }(4\alpha +\mu )s^{3}\gamma ^{3}\varphi
^{3}d_{2}\left( \psi \right) \left\vert z\right\vert ^{2}dxdydt, \\
X_{6} &=&2\func{Im}\int\nolimits_{\Omega }s\gamma \left( 4\alpha +\mu
\right) \partial _{t}z\overline{z}\varphi dxdydt \\
&&-\int\nolimits_{\Omega }\left( 4\alpha +\mu \right) s\gamma ^{2}\varphi
\left\vert z\right\vert ^{2}d_{1}\left( \psi \right) dxdydt \\
&&-\int\nolimits_{\Omega }\left( 4\alpha +\mu \right) s^{2}\gamma
^{3}\varphi \left\vert z\right\vert ^{2}d_{2}\left( \psi \right) dxdydt \\
&&-2\int\nolimits_{\Omega }\left( 4\alpha +\mu \right) s^{2}\gamma
^{3}\varphi ^{2}\left\vert z\right\vert ^{2}d_{2}\left( \psi \right) dxdydt.
\end{eqnarray*}%
By adding (\ref{3.29}) and (\ref{3.31}) we have%
\begin{eqnarray}
&&2\func{Re}\left( P_{s}^{+}z,P_{s}^{-}z\right) _{L^{2}\left( \Omega \right)
}-2\func{Re}\int\nolimits_{\Omega }\left( 4\alpha +\mu \right) s\gamma
\left( P_{s}^{+}z+P_{s}^{-}z\right) \varphi \overline{z}dxdydt  \notag \\
&\geq &2\mu \underset{\Omega }{\int }s\gamma \varphi \left\vert \nabla
_{y}z\right\vert ^{2}dxdydt+\left( 8-8\alpha -2\mu \right) \underset{\Omega }%
{\int }s\gamma \varphi \left\vert \nabla _{x}z\right\vert ^{2}dxdydt  \notag
\\
&&+64\delta _{0}^{2}\int_{\Omega }s^{3}\gamma ^{4}\varphi
^{3}|z|^{2}dxdydt+B_{0}+X_{1}+X_{2}+X_{5}+X_{6}.  \TCItag{3.32}  \label{3.32}
\end{eqnarray}%
On the other hand, since%
\begin{equation*}
\left\Vert P_{s}z+is\partial _{t}\varphi z\right\Vert _{L^{2}\left( \Omega
\right) }^{2}=\left\Vert P_{s}^{+}z\right\Vert _{L^{2}\left( \Omega \right)
}^{2}+\left\Vert P_{s}^{-}z\right\Vert _{L^{2}\left( \Omega \right) }^{2}+2%
\func{Re}\left( P_{s}^{+}z,P_{s}^{-}z\right) _{L^{2}\left( \Omega \right) },
\end{equation*}%
and%
\begin{eqnarray*}
&&-2\func{Re}\int\nolimits_{\Omega }\left( 4\alpha +\mu \right) s\gamma
\left( P_{s}z+is\partial _{t}\varphi z\right) \varphi \overline{z}dxdydt \\
&\leq &\left( 4\alpha +\mu \right) \int_{\Omega }\left\vert
P_{s}z+is\partial _{t}\varphi z\right\vert
^{2}dxdydt+C_{1}s^{2}\int\nolimits_{\Omega }\left\vert z\right\vert
^{2}dxdydt
\end{eqnarray*}%
we have%
\begin{eqnarray*}
C_{2}\int\nolimits_{\Omega }\left\vert P_{s}z+is\partial _{t}\varphi
z\right\vert ^{2}dxdydt &\geq &\int_{\Omega }\left\vert
P_{s}^{+}z\right\vert ^{2}dxdydt+\int_{\Omega }\left\vert
P_{s}^{-}z\right\vert ^{2}dxdydt \\
&&+2s\gamma C_{3}\int_{\Omega }\varphi \left\vert \nabla _{x}z\right\vert
^{2}dxdydt+2s\gamma C_{4}\int_{\Omega }\varphi \left\vert \nabla
_{y}z\right\vert ^{2}dxdydt \\
&&+64\delta _{0}^{2}s^{3}\gamma ^{4}\int\nolimits_{\Omega }\varphi
^{3}\left\vert z\right\vert ^{2}dxdydt+B_{0}+X_{1}+X_{2}+X_{5}+X_{6}.
\end{eqnarray*}%
We see that there exists a constant $\gamma _{1}$ such that for arbitrary $%
\gamma >\gamma _{1},$ the terms of $X_{1}$ and $X_{5}$ can be absorbed by $%
64\delta _{0}^{2}s^{3}\gamma ^{4}\int\nolimits_{\Omega }\varphi
^{3}\left\vert z\right\vert ^{2}dxdydt$, and we have%
\begin{eqnarray*}
C_{2}\int\nolimits_{\Omega }\left\vert P_{s}z+is\partial _{t}\varphi
z\right\vert ^{2}dxdydt &\geq &\int_{\Omega }\left\vert
P_{s}^{+}z\right\vert ^{2}dxdydt+\int_{\Omega }\left\vert
P_{s}^{-}z\right\vert ^{2}dxdydt \\
&&+2s\gamma C_{3}\int_{\Omega }\varphi \left\vert \nabla _{x}z\right\vert
^{2}dxdydt+2s\gamma C_{4}\int_{\Omega }\varphi \left\vert \nabla
_{y}z\right\vert ^{2}dxdydt \\
&&+64\delta _{0}^{2}s^{3}\gamma ^{4}\int\nolimits_{\Omega }\varphi
^{3}\left\vert z\right\vert ^{2}dxdydt+B_{0}+X_{2}+X_{6}.
\end{eqnarray*}%
Since $\varphi >0$ on $\bar{\Omega}$ for $\gamma >\gamma _{1},$ there exist
constants $C_{5}=C_{5}\left( \gamma \right) $ and $s_{1}=s_{1}\left( \gamma
\right) $ such that for all $s>s_{1},$%
\begin{eqnarray*}
C_{2}\int\nolimits_{\Omega }\left\vert P_{s}z+is\partial _{t}\varphi
z\right\vert ^{2}dxdydt &\geq &\int_{\Omega }\left\vert
P_{s}^{+}z\right\vert ^{2}dxdydt+\int_{\Omega }\left\vert
P_{s}^{-}z\right\vert ^{2}dxdydt \\
&&+C_{5}\left( \gamma \right) s\int_{\Omega }\left\vert \nabla
_{x}z\right\vert ^{2}dxdydt+C_{5}\left( \gamma \right) s\int_{\Omega
}\left\vert \nabla _{y}z\right\vert ^{2}dxdydt \\
&&+C_{5}\left( \gamma \right) s^{3}\int\nolimits_{\Omega }\left\vert
z\right\vert ^{2}dxdydt+B_{0}+X_{2}+X_{6}.
\end{eqnarray*}%
Then we choose $s_{2}=s_{2}\left( \gamma \right) >0$ such that $\forall
s>s_{2}$ all the terms of $X_{2}$ and $X_{6}$ can be absorbed into $%
\left\Vert P_{s}^{+}z\right\Vert _{L^{2}\left( \Omega \right) }^{2},$ $%
\left\Vert P_{s}^{-}z\right\Vert _{L^{2}\left( \Omega \right) }^{2},$ $%
C_{5}\left\Vert \nabla _{x}z\right\Vert _{L^{2}\left( \Omega \right) }^{2},$ 
$C_{5}\left\Vert \nabla _{y}z\right\Vert _{L^{2}\left( \Omega \right) }^{2}$
and $C_{5}s^{3}\left\Vert z\right\Vert _{L^{2}\left( \Omega \right) }^{2}.$
Therefore $\int\nolimits_{\Omega }\left\vert P_{s}z+is\partial _{t}\varphi
z\right\vert ^{2}dxdydt\leq 2\int_{\Omega }\left\vert P_{s}z\right\vert
^{2}dxdydt+C_{6}s^{2}\int\nolimits_{\Omega }\left\vert z\right\vert
^{2}dxdydt,$ taking $s>0$ sufficiently large, we have%
\begin{eqnarray*}
C_{7}\int_{\Omega }\left\vert P_{s}z\right\vert ^{2}dxdydt &\geq
&\int_{\Omega }\left\vert P_{s}^{+}z\right\vert ^{2}dxdydt+\int_{\Omega
}\left\vert P_{s}^{-}z\right\vert ^{2}dxdydt \\
&&+s\int_{\Omega }\left\vert \nabla _{x}z\right\vert
^{2}dxdydt+s\int_{\Omega }\left\vert \nabla _{y}z\right\vert ^{2}dxdydt \\
&&+s^{3}\int\nolimits_{\Omega }\left\vert z\right\vert ^{2}dxdydt+B_{0}.
\end{eqnarray*}%
Since $z=0$ on $\Gamma _{x}\cup \Gamma _{y},$ $\nabla _{y}z=0$ and $\nabla
_{x}z=\left( \partial _{\nu }z\right) \cdot \nu $ on $\Gamma _{x},$ all the
integrations on $\Gamma _{y}$ vanish in (\ref{3.28}), we have%
\begin{eqnarray*}
B_{0} &=&-4\func{Re}\int_{\Gamma _{x}}s\gamma \varphi \left( \partial _{\nu
}z\right) \nabla _{x}\psi \cdot \nabla _{x}\overline{z}dS_{x}dydt+\int%
\nolimits_{\Gamma _{x}}2s\gamma \varphi \left( \partial _{\nu }\psi \right)
\left\vert \nabla _{x}z\right\vert ^{2}dS_{x}dydt \\
&=&-8\int_{\Gamma _{x}}s\gamma \varphi \left\vert \partial _{\nu
}z\right\vert ^{2}(x-x_{0})\cdot \nu dS_{x}dydt+4\int\nolimits_{\Gamma
_{x}}s\gamma \varphi \left\vert \partial _{\nu }z\right\vert
^{2}(x-x_{0})\cdot \nu dS_{x}dydt \\
&=&-4\int\nolimits_{\Gamma _{x}}s\gamma \varphi \left\vert \partial _{\nu
}z\right\vert ^{2}(x-x_{0})\cdot \nu dS_{x}dydt \\
&\geq &-4\int\nolimits_{\Gamma _{x}\cap \left\{ (x-x_{0})\cdot \nu \geq
0\right\} }s\gamma \varphi \left\vert \partial _{\nu }z\right\vert
^{2}(x-x_{0})\cdot \nu dS_{x}dydt.
\end{eqnarray*}%
We obtain%
\begin{eqnarray*}
&&\int_{\Omega }\left\vert P_{s}^{+}z\right\vert ^{2}dxdydt+\int_{\Omega
}\left\vert P_{s}^{-}z\right\vert ^{2}dxdydt+s\int_{\Omega }\left\vert
\nabla _{x}z\right\vert ^{2}dxdydt \\
&&+s\int_{\Omega }\left\vert \nabla _{y}z\right\vert
^{2}dxdydt+s^{3}\int\nolimits_{\Omega }\left\vert z\right\vert ^{2}dxdydt. \\
&\leq &C_{8}\int_{\Omega }\left\vert P_{s}z\right\vert
^{2}dxdydt+C_{8}s\int_{\partial D_{+}\times G\times (-T,T)}\left\vert
\partial _{\nu }z\right\vert ^{2}dS_{x}dydt.
\end{eqnarray*}

Finally, we rewrite our inequality with $z$ instead of $u$. By the relation%
\begin{eqnarray*}
\left\vert z\right\vert ^{2} &=&e^{2s\varphi }\left\vert u\right\vert ^{2},%
\text{ }\left\vert \partial _{\nu }z\right\vert ^{2}=\left\vert \partial
_{\nu }u\right\vert ^{2}e^{2s\varphi }\text{ on }\partial D_{+}\times
G\times (-T,T), \\
\left\vert \nabla _{x}ue^{s\varphi }\right\vert ^{2} &=&\left\vert \nabla
_{x}z-s\lambda \varphi e^{s\varphi }u\nabla _{x}\psi \right\vert ^{2}\leq
2\left\vert \nabla _{x}z\right\vert ^{2}+2s^{2}\lambda ^{2}\varphi
^{2}\left\vert \nabla _{x}\psi \right\vert ^{2}\left\vert z\right\vert ^{2},
\\
\left\vert \nabla _{y}ue^{s\varphi }\right\vert ^{2} &=&\left\vert \nabla
_{y}z-s\lambda \varphi e^{s\varphi }u\nabla _{y}\psi \right\vert ^{2}\leq
2\left\vert \nabla _{y}z\right\vert ^{2}+2s^{2}\lambda ^{2}\varphi
^{2}\left\vert \nabla _{y}\psi \right\vert ^{2}\left\vert z\right\vert ^{2},
\\
\left\vert L_{0}u\right\vert ^{2} &\leq &2\left\vert Lu\right\vert
^{2}+2\left\vert \sum\limits_{i=1}^{n}a_{i}\left( x,y,t\right)
u_{x_{i}}+\sum\limits_{j=1}^{m}b_{j}\left( x,y,t\right)
u_{y_{j}}+a_{0}\left( x,y,t\right) u\left( x,y,t\right) \right\vert ^{2},
\end{eqnarray*}%
we see that there exist positive constants $C_{9}=C_{9}\left( \gamma \right) 
$ and $s_{0}>s_{2}\left( \gamma \right) $ such that for all $s>s_{0},$%
\begin{eqnarray*}
&&\int_{\Omega }(s\left\vert \nabla _{y}u\right\vert ^{2}+s\left\vert \nabla
_{x}u\right\vert ^{2}+s^{3}\left\vert u\right\vert ^{2})e^{2s\varphi }dxdydt
\\
&\leq &C_{9}\int_{\Omega }\left\vert Lu\right\vert ^{2}e^{2s\varphi
}dxdydt+C_{9}\int_{\partial D_{+}\times G\times (-T,T)}s\left\vert \partial
_{\nu }u\right\vert ^{2}e^{2s\varphi }dS_{x}dydt.
\end{eqnarray*}

We complete the proof of Theorem 1.

\section{\textbf{Proof of Theorem 1}}

Since $u$ itself does not satisfy (\ref{3.6}), in order to apply Proposition
1, we have to introduce a cut-off function. Moreover, we need to introduce
several notations. We set%
\begin{eqnarray}
\widetilde{r} &=&\max_{x\in \overline{D}}\left\vert x-x_{0}\right\vert , 
\TCItag{4.1}  \label{4.1} \\
r &=&\min_{x\in \overline{D}}\left\vert x-x_{0}\right\vert .  \TCItag{4.2}
\label{4.2}
\end{eqnarray}%
By $x_{0}\notin \overline{D}$, we see that $r>0$. We choose $\rho >1$
sufficiently large so that%
\begin{equation}
\frac{\widetilde{r}}{r}<\rho .  \tag{4.3}  \label{4.3}
\end{equation}%
By (\ref{4.3}) and the assumption on $L,$ we have%
\begin{equation}
\frac{\alpha L^{2}}{\rho ^{2}}<r^{2}<\widetilde{r}^{2}<\alpha L^{2}. 
\tag{4.4}  \label{4.4}
\end{equation}%
Furthermore, if necessary, we choose smaller $\alpha $ such that%
\begin{equation}
r^{2}>\alpha ^{2}L^{2}.  \tag{4.5}  \label{4.5}
\end{equation}%
We arbitrarily choose $y_{0}=\left( y_{1},y_{2},\ldots ,y_{m}\right) \in 
\mathbb{R}
^{m}$ satisfying%
\begin{equation}
\left\vert y_{0}\right\vert \leq L-\frac{L}{\rho }-\epsilon .  \tag{4.6}
\label{4.6}
\end{equation}%
We set%
\begin{eqnarray*}
G_{1} &=&\left\{ y\in 
\mathbb{R}
^{m};\text{ }\left\vert y-y_{0}\right\vert <L\right\} , \\
G_{2} &=&\left\{ y\in 
\mathbb{R}
^{m};\text{ }\left\vert y\right\vert <2L\right\} , \\
\Omega _{0} &=&D\times G_{1}\times \left\{ t=0\right\} ,\text{ }\Omega
_{1}=D\times G_{1}\times (-T,T), \\
\Omega _{2} &=&D\times G_{2}\times (-T,T).
\end{eqnarray*}%
Then (\ref{4.1}), (\ref{4.2}) and (\ref{4.4}) yields, if $x\in D$ and $%
\left\vert y-y_{0}\right\vert \leq L,$%
\begin{eqnarray}
\psi \left( x,y,\mp T\right) &=&\left\vert x-x_{0}\right\vert ^{2}-\alpha
\left\vert y-y_{0}\right\vert ^{2}-\beta T^{2}  \notag \\
&\leq &\left\vert x-x_{0}\right\vert ^{2}-\alpha \left\vert
y-y_{0}\right\vert ^{2}  \notag \\
&\leq &\widetilde{r}^{2}-\alpha L^{2}  \notag \\
&<&0  \TCItag{4.7}  \label{4.7}
\end{eqnarray}%
and if $x\in D,$ $\left\vert y-y_{0}\right\vert \leq L$ and $\left\vert
t\right\vert <T,$%
\begin{eqnarray}
\psi \left( x,y,t\right) &=&\left\vert x-x_{0}\right\vert ^{2}-\alpha
\left\vert y-y_{0}\right\vert ^{2}-\beta \left\vert t\right\vert ^{2}  \notag
\\
&\leq &\left\vert x-x_{0}\right\vert ^{2}-\alpha \left\vert
y-y_{0}\right\vert ^{2}  \notag \\
&\leq &\widetilde{r}^{2}-\alpha L^{2}  \notag \\
&<&0  \TCItag{4.8}  \label{4.8}
\end{eqnarray}%
and if $x\in D$ and $\left\vert y-y_{0}\right\vert \leq \frac{L}{\rho },$%
\begin{eqnarray}
\psi \left( x,y,0\right) &=&\left\vert x-x_{0}\right\vert ^{2}-\alpha
\left\vert y-y_{0}\right\vert ^{2}  \notag \\
&\geq &r^{2}-\alpha \frac{L^{2}}{p^{2}}  \notag \\
&>&0.  \TCItag{4.9}  \label{4.9}
\end{eqnarray}%
Therefore, for small $\epsilon >0$ there exists $\delta >0$ such that 
\begin{equation}
\psi \left( x,y,t\right) <-\epsilon ,\text{ }x\in D,  \tag{4.10}
\label{4.10}
\end{equation}%
if $T-2\delta \leq \left\vert t\right\vert \leq T$ or $L-2\delta \leq
\left\vert y-y_{0}\right\vert \leq L$ and 
\begin{equation}
\psi \left( x,y,t\right) >\epsilon ,\text{ }x\in D,\text{ }\left\vert
t\right\vert <\delta ,\text{ }\left\vert y-y_{0}\right\vert \leq \frac{L}{%
\rho }.  \tag{4.11}  \label{4.11}
\end{equation}%
Let us define a cut-off function $\chi \left( y,t\right) =\chi _{0}\left(
t\right) \chi _{0}\left( \left\vert y-y_{0}\right\vert \right) ,$ where $%
\chi _{0}\in C_{0}^{\infty }\left( 
\mathbb{R}
\right) $ such that $0\leq \chi _{0}\leq 1$ and%
\begin{equation*}
\chi _{0}\left( t\right) =\left\{ 
\begin{array}{c}
0, \\ 
1,%
\end{array}%
\begin{array}{c}
T-\delta \leq \left\vert t\right\vert \leq T, \\ 
\left\vert t\right\vert \leq T-2\delta ,%
\end{array}%
\right.
\end{equation*}%
\begin{equation*}
\chi _{0}\left( \left\vert y-y_{0}\right\vert \right) =\left\{ 
\begin{array}{c}
0, \\ 
1,%
\end{array}%
\begin{array}{c}
L-\delta \leq \left\vert y-y_{0}\right\vert \leq L, \\ 
\left\vert y-y_{0}\right\vert \leq L-2\delta .%
\end{array}%
\right.
\end{equation*}%
Then we see that $\chi \in C_{0}^{\infty }\left( 
\mathbb{R}
^{m+1}\right) ,$ $0\leq \chi \leq 1$ and%
\begin{equation}
\chi \left( y,t\right) =\left\{ 
\begin{array}{c}
0, \\ 
1,%
\end{array}%
\begin{array}{c}
T-\delta \leq \left\vert t\right\vert \leq T\text{ or }L-\delta \leq
\left\vert y-y_{0}\right\vert \leq L, \\ 
\left\vert t\right\vert \leq T-2\delta \text{ and }\left\vert
y-y_{0}\right\vert \leq L-2\delta .%
\end{array}%
\right.  \tag{4.12}  \label{4.12}
\end{equation}%
By choosing $\delta >0$ smaller if necessary, we assume%
\begin{equation}
\frac{L}{\rho }<L-2\delta .  \tag{4.13}  \label{4.13}
\end{equation}%
\qquad We set%
\begin{equation*}
w_{k}=\left( \partial _{t}^{k}u\right) \chi ,\text{ }k=1,2.
\end{equation*}%
Then%
\begin{eqnarray}
Aw_{k} &=&f\partial _{t}^{k}R\chi +2\left( \nabla _{y}\partial
_{t}^{k}u\cdot \nabla _{y}\chi \right) +\partial _{t}^{k}u\left( \Delta
_{y}\chi +i\partial _{t}\chi \right) ,  \notag \\
x &\in &D,\text{ }y\in G_{1},\text{ }k=1,2  \TCItag{4.14}  \label{4.14}
\end{eqnarray}%
and%
\begin{eqnarray}
w_{k}\left( x,y,t\right) &=&\left\vert \nabla _{y}w_{k}\left( x,y,t\right)
\right\vert =0,\text{ }\left( x,y,t\right) \in D\times \partial G_{1}\times
(-T,T),  \notag \\
w_{k}\left( x,y,t\right) &=&0,\text{ }\left( x,y,t\right) \in \partial
D\times G_{1}\times (-T,T),  \notag \\
w_{k}\left( x,y,T\right) &=&w_{k}\left( x,y,-T\right) =0,\text{ }\left(
x,y,t\right) \in D\times G_{1}.  \TCItag{4.15}  \label{4.15}
\end{eqnarray}%
From (\ref{4.6}) we note that%
\begin{equation}
G_{1}\subset G_{2}.  \tag{4.16}  \label{4.16}
\end{equation}%
By (\ref{4.14})-(\ref{4.15}), we can apply the Carleman estimate (see
Proposition 1) to $w_{1},w_{2}:$%
\begin{eqnarray}
&&\int\nolimits_{\Omega _{1}}\sum\limits_{k=1}^{2}\left( s\left\vert \nabla
_{x}w_{k}\right\vert ^{2}+s\left\vert \nabla _{y}w_{k}\right\vert
^{2}+s\left\vert w_{k}\right\vert ^{2}\right) e^{2s\varphi }dxdydt  \notag \\
&\leq &C\int\nolimits_{\Omega _{1}}\sum\limits_{k=1}^{2}\chi
^{2}f^{2}\left\vert \partial _{t}^{k}R\right\vert ^{2}e^{2s\varphi }dxdydt 
\notag \\
&&+C\int\nolimits_{\Omega _{1}}\sum\limits_{k=1}^{2}\left( \left\vert
2\left( \nabla _{y}\partial _{t}^{k}u\cdot \nabla _{y}\chi \right) +\partial
_{t}^{k}u\left( \Delta _{y}\chi +i\partial _{t}\chi \right) \right\vert
^{2}\right) e^{2s\varphi }dxdydt  \notag \\
&&+C\int\nolimits_{\partial D_{+}\times G_{1}\times \left( -T,T\right)
}\sum\limits_{k=1}^{2}s\left\vert \partial _{\nu }w_{k}\right\vert
^{2}e^{2s\varphi }dS_{x}dydt  \notag \\
&=&S_{1}+S_{2}+S_{3}.  \TCItag{4.17}  \label{4.17}
\end{eqnarray}

Here and henceforth, $C>0$ denotes a generic constant which is independent
of $s>0$. From the assumption on $R$, we have%
\begin{eqnarray*}
S_{1} &=&C\int\limits_{\Omega _{1}}\sum\limits_{k=1}^{2}\chi
^{2}f^{2}\left\vert \partial _{t}^{k}R\right\vert ^{2}e^{2s\varphi }dxdydt \\
&\leq &C\int\limits_{\Omega _{1}}\chi ^{2}f^{2}e^{2s\varphi }dxdydt.
\end{eqnarray*}%
By (\ref{4.12}), we see that $\partial _{t}\chi =0$ for $\left\vert
t\right\vert \leq T-2\delta $ or $T-\delta \leq \left\vert t\right\vert \leq
T$ and $\left\vert \nabla _{y}\chi \right\vert =\Delta _{y}\chi =0$ for $%
\left\vert y-y_{0}\right\vert \leq L-2\delta $ or $L-\delta \leq \left\vert
y-y_{0}\right\vert \leq L.$ Therefore, if $\left\vert t\right\vert \in \left[
0,T-2\delta \right] \cup \left[ T-\delta ,T\right] $ and $\left\vert
y-y_{0}\right\vert \in \left[ 0,L-2\delta \right] \cup \left[ L-\delta ,L%
\right] $, then $\partial _{t}\chi =\left\vert \nabla _{y}\chi \right\vert
=\Delta _{y}\chi =0$. Hence%
\begin{eqnarray}
S_{2} &=&C\left( \int\nolimits_{\left\{ T-2\delta \leq \left\vert
t\right\vert \leq T-\delta \right\} \cap \Omega
_{1}}\sum\limits_{k=1}^{2}\left\vert 2\left( \nabla _{y}\partial
_{t}^{k}u\cdot \nabla _{y}\chi \right) \right. \right.  \notag \\
&&\left. \left. +\partial _{t}^{k}u\left( \Delta _{y}\chi +i\partial
_{t}\chi \right) \right\vert ^{2}e^{2s\varphi }dxdydt\right)  \notag \\
&&+C\left( \int\nolimits_{\left\{ L-2\delta \leq \left\vert
y-y_{0}\right\vert \leq L-\delta \right\} \cap \Omega
_{1}}\sum\limits_{k=1}^{2}\left\vert 2\left( \nabla _{y}\partial
_{t}^{k}u\cdot \nabla _{y}\chi \right) \right. \right.  \notag \\
&&\left. \left. +\partial _{t}^{k}u\left( \Delta _{y}\chi +i\partial
_{t}\chi \right) \right\vert ^{2}e^{2s\varphi }dxdydt\right)  \notag \\
&\mathbf{=}&C\left( \int\nolimits_{\left\{ T-2\delta \leq \left\vert
t\right\vert \leq T-\delta \right\} \cap \Omega
_{1}}\sum\limits_{k=1}^{2}\left\vert 2\left( \nabla _{y}\partial
_{t}^{k}u\cdot \nabla _{y}\chi \right) \right. \right.  \notag \\
&&\left. \left. +\partial _{t}^{k}u\left( \Delta _{y}\chi +i\partial
_{t}\chi \right) \right\vert ^{2}\left( \exp (2se^{-\gamma \epsilon
})\right) )dxdydt\right)  \notag \\
&&\mathbf{+}C\left( \int\nolimits_{\left\{ L-2\delta \leq \left\vert
y-y_{0}\right\vert \leq L-\delta \right\} \cap \Omega
_{1}}\sum\limits_{k=1}^{2}\left\vert 2\left( \nabla _{y}\partial
_{t}^{k}u\cdot \nabla _{y}\chi \right) \right. \right.  \notag \\
&&\left. \left. +\partial _{t}^{k}u\left( \Delta _{y}\chi +i\partial
_{t}\chi \right) \right\vert ^{2}\left( \exp (2se^{-\gamma \epsilon
})\right) )dxdydt\right)  \notag \\
&\leq &C\int\nolimits_{\Omega _{1}}\sum\limits_{k=1}^{2}\left( \left\vert
\partial _{t}^{k}u\right\vert ^{2}+\left\vert \nabla _{y}\partial
_{t}^{k}u\right\vert ^{2}\right) e^{2s\kappa _{1}}dxdydt  \notag \\
&\leq &CM^{2}e^{2s\kappa _{1}}.  \TCItag{4.18}  \label{4.18}
\end{eqnarray}%
Here and henceforth we set%
\begin{equation*}
\kappa _{1}=e^{-\gamma \epsilon },\text{ }\kappa _{2}=e^{\gamma \epsilon }.
\end{equation*}%
Finally, we have%
\begin{eqnarray}
S_{3} &=&C\int\nolimits_{\partial D_{+}\times G_{1}\times \left( -T,T\right)
}\sum\limits_{k=1}^{2}s\left\vert \partial _{\nu }w_{k}\right\vert
^{2}e^{2s\varphi }dS_{x}dydt  \notag \\
&\leq &Ce^{cs}\int\nolimits_{\partial D_{+}\times G_{2}\times \left(
-T,T\right) }\sum\limits_{k=1}^{2}\left\vert \partial _{\nu }\partial
_{t}^{k}u\right\vert ^{2}dS_{x}dydt:=Ce^{cs}d^{2},  \TCItag{4.19}
\label{4.19}
\end{eqnarray}%
where%
\begin{equation*}
d^{2}=\int\nolimits_{\partial D_{+}\times G_{2}\times \left( -T,T\right)
}\sum\limits_{k=1}^{2}\left\vert \partial _{\nu }\partial
_{t}^{k}u\right\vert ^{2}dS_{x}dydt.
\end{equation*}%
As a result, (\ref{4.17}) yields%
\begin{eqnarray}
&&\int\limits_{\Omega _{1}}\sum\limits_{k=1}^{2}\left( s\left\vert \nabla
_{x}w_{k}\right\vert ^{2}+s\left\vert \nabla _{y}w_{k}\right\vert
^{2}+s^{3}\left\vert w_{k}\right\vert ^{2}\right) e^{2s\varphi }dxdydt 
\notag \\
&\leq &C\int\limits_{\Omega _{1}}\chi ^{2}f^{2}e^{2s\varphi
}dxdydt+CM^{2}e^{2s\kappa _{1}}+Ce^{cs}d^{2}.  \TCItag{4.20}  \label{4.20}
\end{eqnarray}%
Now, by using the fact that $\chi (y,-T)=0$ for $y\in G_{1}$ by (\ref{4.12}%
), we can write%
\begin{eqnarray*}
&&\int\nolimits_{\Omega _{0}}\left\vert \chi \left( y,0\right) \right\vert
^{2}\left\vert i\partial _{t}u\left( x,y,0\right) \right\vert
^{2}e^{2s\varphi \left( x,y,0\right) }dxdy \\
&=&\int\limits_{-T}^{0}\partial _{t}\left( \int\limits_{\Omega _{0}}\chi
^{2}\left\vert \partial _{t}u\left( x,y,t\right) \right\vert
^{2}e^{2s\varphi \left( x,y,t\right) }dxdy\right) dt \\
&=&\int\limits_{-T}^{0}\int\nolimits_{\Omega _{0}}\left( 2\chi \partial
_{t}\chi \left\vert \partial _{t}u\right\vert ^{2}+\chi ^{2}\partial
_{t}\left( \left\vert \partial _{t}u\right\vert ^{2}\right) +\chi
^{2}\left\vert \partial _{t}u\right\vert ^{2}2s\partial _{t}\varphi \right)
e^{2s\varphi \left( x,y,t\right) }dxdydt \\
&=&\int\limits_{-T}^{0}\int\nolimits_{\Omega _{0}}\left( 2\chi \partial
_{t}\chi \left\vert \partial _{t}u\right\vert ^{2}+2\chi ^{2}\func{Re}%
\partial _{t}^{2}u\partial _{t}\bar{u}+\chi ^{2}\left\vert \partial
_{t}u\right\vert ^{2}2s\partial _{t}\varphi \right) e^{2s\varphi \left(
x,y,t\right) }dxdydt \\
&\leq &\int\limits_{-T}^{T}\int\nolimits_{\Omega _{0}}\left\vert \chi
\partial _{t}u\right\vert ^{2}e^{2s\varphi
}dxdydt+\int\limits_{-T}^{T}\int\nolimits_{D\times G_{1}}\left\vert \partial
_{t}\chi \partial _{t}u\right\vert ^{2}e^{2s\varphi }dxdydt \\
&&+\int\limits_{-T}^{T}\int\limits_{D\times G_{1}}\left\vert \chi \partial
_{t}u\right\vert ^{2}e^{2s\varphi
}dxdydt+\int\limits_{-T}^{T}\int\limits_{\Omega _{0}}\left\vert \chi
\partial _{t}^{2}u\right\vert ^{2}e^{2s\varphi }dxdydt \\
&&+s\int\limits_{-T}^{T}\int\limits_{\Omega _{0}}\left\vert \chi \partial
_{t}u\right\vert ^{2}e^{2s\varphi }dxdydt
\end{eqnarray*}%
\begin{eqnarray*}
&\leq &\int\limits_{-T}^{T}\int\limits_{\Omega _{0}}\left\vert \partial
_{t}\chi \partial _{t}u\right\vert ^{2}e^{2s\varphi }dxdydt \\
&&+C\int\limits_{-T}^{T}\int\limits_{\Omega _{0}}\left( \left\vert \chi
\partial _{t}u\right\vert ^{2}+\left\vert \chi \partial _{t}^{2}u\right\vert
^{2}+s\left\vert \chi \partial _{t}u\right\vert ^{2}\right) e^{2s\varphi
}dxdydt \\
&\leq &CM^{2}e^{2s\kappa _{1}}+C\int\limits_{\Omega _{1}}\left( \left\vert
w_{1}\right\vert ^{2}+\left\vert w_{2}\right\vert ^{2}+s\left\vert
w_{1}\right\vert ^{2}\right) e^{2s\varphi }dxdydt.
\end{eqnarray*}%
Hence%
\begin{eqnarray*}
&&\int\nolimits_{\Omega _{0}}\left\vert \chi \left( y,0\right) \right\vert
^{2}\left\vert i\partial _{t}u\left( x,y,0\right) \right\vert
^{2}e^{2s\varphi \left( x,y,0\right) }dxdy \\
&\leq &CM^{2}e^{2s\kappa _{1}}+C\int\nolimits_{\Omega _{1}}\left(
s\left\vert w_{1}\right\vert ^{2}+\left\vert w_{2}\right\vert ^{2}\right)
e^{2s\varphi }dxdydt.
\end{eqnarray*}%
Then, applying (\ref{4.20}), we obtain%
\begin{eqnarray}
&&\int\nolimits_{\Omega _{0}}\left\vert \chi _{0}\left( \left\vert
y-y_{0}\right\vert \right) \right\vert ^{2}\left\vert i\partial _{t}u\left(
x,y,0\right) \right\vert ^{2}e^{2s\varphi \left( x,y,0\right) }dxdy  \notag
\\
&\leq &\frac{C}{s}\int\nolimits_{\Omega _{1}}\chi ^{2}\left( y,t\right)
\left\vert f\right\vert ^{2}e^{2s\varphi \left( x,y,t\right)
}dxdydt+CM^{2}e^{2s\kappa _{1}}+Ce^{cs}d^{2}  \notag \\
&=&\frac{C}{s}\int\nolimits_{\Omega _{1}}\chi _{0}^{2}\left( t\right) \chi
_{0}^{2}\left( \left\vert y-y_{0}\right\vert \right) \left\vert f\right\vert
^{2}e^{2s\varphi \left( x,y,t\right) }dxdydt+CM^{2}e^{2s\kappa
_{1}}+Ce^{cs}d^{2}  \notag \\
&\leq &\frac{C}{s}\int\nolimits_{\Omega _{0}}\chi _{0}^{2}\left( \left\vert
y-y_{0}\right\vert \right) \left\vert f\right\vert ^{2}e^{2s\varphi \left(
x,y,0\right) }dxdy+CM^{2}e^{2s\kappa _{1}}+Ce^{cs}d^{2},  \TCItag{4.21}
\label{4.21}
\end{eqnarray}%
where we used $\left\vert \chi _{0}\left( t\right) \right\vert \leq 1$ and $%
e^{2s\varphi \left( x,y,t\right) }\leq e^{2s\varphi \left( x,y,0\right) }$
for $x\in D$ and $y\in G_{1}.$

On the other hand, by substituting $t=0$ in (\ref{2.1}) and applying $%
u\left( x,y,0\right) =0$ and $R\left( x,y,0\right) \neq 0,$ for $x\in 
\overline{D}$ and $\left\vert y\right\vert \leq 2L,$ we get%
\begin{equation}
f\left( x,y\right) =\dfrac{i\partial _{t}u\left( x,y,0\right) }{R\left(
x,y,0\right) },\text{ }x\in \overline{D},\text{ }\left\vert y\right\vert
\leq 2L.  \tag{4.22}  \label{4.22}
\end{equation}%
By applying (\ref{4.22}) in (\ref{4.21}), we have%
\begin{eqnarray*}
&&\int\nolimits_{\Omega _{0}}\chi _{0}^{2}\left( \left\vert
y-y_{0}\right\vert \right) \left\vert f\right\vert ^{2}e^{2s\varphi \left(
x,y,0\right) }dxdy \\
&\leq &\frac{C}{s}\int\nolimits_{\Omega _{0}}\chi _{0}^{2}\left( \left\vert
y-y_{0}\right\vert \right) \left\vert f\right\vert ^{2}e^{2s\varphi \left(
x,y,0\right) }dxdy+CM^{2}e^{2s\kappa _{1}}+Ce^{cs}d^{2}
\end{eqnarray*}%
for all large $s>0$. Now, we absorb the first term on the right-hand side
into the left-hand side by choosing $s>0$ large, we get%
\begin{equation*}
\int\nolimits_{\Omega _{0}}\chi _{0}^{2}\left( \left\vert y-y_{0}\right\vert
\right) \left\vert f\right\vert ^{2}e^{2s\varphi \left( x,y,0\right)
}dxdy\leq CM^{2}e^{2s\kappa _{1}}+Ce^{cs}d^{2}
\end{equation*}%
for all large $s>0$.

Replacing the integration domain on the left-hand side by $D\times \left\{ y;%
\text{ }\left\vert y-y_{0}\right\vert <\dfrac{L}{\rho }\right\} \subset
\Omega _{0}$ and using the facts that $\chi _{0}\left( \left\vert
y-y_{0}\right\vert \right) =1$ in $D\times \left\{ y;\text{ }\left\vert
y-y_{0}\right\vert <\dfrac{L}{\rho }\right\} $ and%
\begin{equation*}
e^{2s\varphi \left( x,y,0\right) }=\exp \left( 2se^{\gamma \psi \left(
x,y,0\right) }\right) >\exp \left( 2se^{\gamma \epsilon }\right)
=e^{2s\kappa _{2}},
\end{equation*}%
we obtain,%
\begin{equation*}
e^{2s\kappa _{2}}\int\nolimits_{D\times \left\{ y;\text{ }\left\vert
y-y_{0}\right\vert <\frac{L}{\rho }\right\} }\left\vert f\right\vert
^{2}dxdy\leq CM^{2}e^{2s\kappa _{1}}+Ce^{cs}d^{2}
\end{equation*}%
for all $s\geq s_{0}$, where $s_{0}$ is some constant. Since $\kappa
_{2}>\kappa _{1},$ the last inequality implies%
\begin{equation}
\int\nolimits_{D\times \left\{ y;\text{ }\left\vert y-y_{0}\right\vert <%
\frac{L}{\rho }\right\} }\left\vert f\right\vert ^{2}dxdy\leq
CM^{2}e^{-2s\kappa }+Ce^{cs}d^{2}  \tag{4.23}  \label{4.23}
\end{equation}%
for all $s\geq s_{0}$, where $\kappa =\kappa _{2}-\kappa _{1}>0$. We
separately consider the two cases:

Case 1. Let $M\geq d.$ Choosing $s\geq 0$ such that

\begin{equation*}
M^{2}e^{-2s\kappa }=e^{Cs}d^{2},\text{ that is},\text{ }s=\frac{2}{C+2\kappa 
}\log \frac{M}{d}\geq 0,
\end{equation*}%
we obtain%
\begin{equation*}
\int\nolimits_{D\times \left\{ y;\text{ }\left\vert y-y_{0}\right\vert <%
\frac{L}{\rho }\right\} }\left\vert f\right\vert ^{2}dxdy\leq 2M^{\frac{2C}{%
C+2\kappa }}d^{\frac{4\kappa }{C+2\kappa }}.
\end{equation*}%
\ \ \ \ \ \ \ 

Case 2. Let $M<d.$ Then setting $s=0$ in (\ref{4.23}) we have%
\begin{equation*}
\int\nolimits_{D\times \left\{ y;\text{ }\left\vert y-y_{0}\right\vert <%
\frac{L}{\rho }\right\} }\left\vert f\right\vert ^{2}dxdy\leq 2Cd^{2}.
\end{equation*}%
Therefore we can choose $\theta \in (0,1)$ such that%
\begin{equation*}
\int\nolimits_{D\times \left\{ y;\text{ }\left\vert y-y_{0}\right\vert <%
\frac{L}{\rho }\right\} }\left\vert f\right\vert ^{2}dxdy\leq C\left(
d^{2\theta }+d^{2}\right)
\end{equation*}%
for all $y_{0}\in 
\mathbb{R}
^{m}$ satisfying $\left\vert y_{0}\right\vert <L-\frac{L}{\rho }-\epsilon .$
By $\left\Vert \partial _{t}u\right\Vert _{H^{2}\left( D\times
\{|y|<2L\}\times (-T,T)\right) }\leq M,$ the trace theorem yields $d\leq CM,$
which implies $d\leq Cd^{\theta }.$ Varying $y_{0}$ and noting%
\begin{equation*}
\bigcup \left\{ y\in 
\mathbb{R}
^{m};\text{ }\left\vert y-y_{0}\right\vert \leq \frac{L}{\rho },\text{ }%
\left\vert y_{0}\right\vert <L-\frac{L}{\rho }-\varepsilon \right\} =\left\{
y\in 
\mathbb{R}
^{m};\text{ }\left\vert y\right\vert <L-\epsilon \right\} ,
\end{equation*}%
we obtain%
\begin{equation*}
\int\nolimits_{D\times \left\{ y;\text{ }\left\vert y\right\vert <L-\epsilon
\right\} }\left\vert f\left( x,y\right) \right\vert ^{2}dxdy\leq Cd^{2\theta
}.
\end{equation*}%
Thus the proof of Theorem 1 is completed.

\end{document}